\magnification=\magstep1
\def\A{{\cal{A}}}
\def\B{{\cal{B}}}
\def\cc{{\cal{C}}}
\def\E{{\cal{E}}}
\def\F{{\cal{F}}}
\def\G{{\cal{G}}}
\def\H{{\cal{H}}}
\def\I{{\cal I}}
\def\J{{\cal J}}
\def\K{{\cal K}}
\def\L{{\cal L}}
\def\P{{~p}}
\def\R{{\cal R}}
\def\U{{\cal U}}
\def\V{{\cal V}}
\def\Y{{\cal Y}}
\def\calo{{\cal O}}
\def\wt{\widetilde}
\def\tb{\widetilde{\B}}
\def\tc{\widetilde{c}}
\def\te{\widetilde{\E}}
\def\tg{\widetilde{G}}
\def\tn{\widetilde{N}}
\def\tq{\widetilde{Q}}
\def\tr{\widetilde{\R}}
\def\tt{\widetilde{T}}
\def\txi{\widetilde\xi}
\def\tz{\widetilde z}
\def\wth{\widetilde h}
\def\teta{\widetilde\eta}
\def\tth{\widetilde\theta}
\def\tnu{\widetilde\nu}
\def\bsl{\backslash}
\def\sprod{\T}
\def\vp{\varphi}
\def\vt{\vartheta}
\def\rank{\rm{rank}}

\def\limsup{\mathop{\overline{\rm lim}}}


\mathsurround=1.5pt
\parskip=5pt

\font\math=msbm10
\font\submath=msbm10 at 7pt
\def\bo#1{\hbox{\math#1}}
\def\bs#1{\hbox{\submath#1}}

\font\bigtenrm=cmr10 scaled\magstep3
\font\medtenrm=cmr10 scaled\magstep1

\def\det{{\rm det}\,}
\def\al{\alpha}
\def\be{\beta}
\def\gam{\gamma}
\def\Gam{\Gamma}

\def\de{\delta}
\def\del{\delta}

\def\lam{\lambda}
\def\Om{\Omega}
\def\rhob{\rho_{\beta}}

\def\sig{\sigma}

\def\Proof:{\noindent{\bf Proof:}}

\def\r{\bo R}
\def\C{\bo C}
\def\N{\bo N}
\def\T{\bo T}
\def\Z{\bo Z}
\def\({\left(}
\def\){\right)}
\def\ovl{\overline}
\def\oR{\ovl{\R}}
\def\z{\ovl{z}}
\def\x{\ovl{x}}
\def\oc{\ovl{c}}
\def\oy{\ovl{y}}
\def\oY{\ovl{\cal Y}}
\def\oxi{\ovl{\xi}}
\def\oeta{\ovl{\eta}}

\def\spo{Sp_0}
\def\sph{{Sp(h)}}
\def\sphd{Sp(h_\gam(\del))}
\def\spoh{Sp_0(h_{\gam^{-1}}(\del))}

\def\hgodl{h_{\gam^{-1}}(\del)}
\def\hgd{h_{\gam}(\del)}


\def\qed{{\unskip\nobreak\hfil\penalty50\hskip .001pt
          \hbox{ }\nobreak\hfil
          \vrule height 1.2ex width 1.1ex depth -.1ex
          \parfillskip=0pt\finalhyphendemerits=0\medbreak}\rm}

\def\qeddis{\eqno{\vrule height 1.2ex width 1.1ex depth -.1ex}
                   $$\medbreak\rm}

\def\om{\omega}
\def\h{\widehat{h}}

\def\ovl{\overline}



\centerline{\bf \bigtenrm Exponentially decaying eigenvectors for certain}

\centerline{\bf \bigtenrm almost periodic operators}

\bigskip
\bigskip

\centerline{\medtenrm \bf Norbert Riedel}

\bigskip
\bigskip

{\narrower\medskip\noindent {\bf Abstract.} For every point $\chi$ in the
spectrum of the operator
$$(h(\delta)\xi)=\xi_{n+1}+\xi_{n-1}+\beta\left(\delta e^{2\pi\al ni}
+\delta^{-1}e^{-2\pi\al ni}\right)\xi_n$$ on $\ell^2(\Z)$ there exists a
complex number $x$ of modulus one such that the equation
$$\xi_{n+1}+\xi_{n-1}+\beta\left(x\delta e^{2\pi\al ni}
+\x\delta^{-1}e^{-2\pi\al ni}\right)\xi_n=\chi\xi_n$$ has a non-trivial
solution satisfying the condition
$$\limsup_{|n|\to\infty}|\xi_n|^{1/|n|}\le\delta^{-1}\beta^{-1}$$
provided that $\beta,\delta>1$ and $\al$ satisfies the diophantine
condition
$$\lim_{n\to\infty}|\sin\pi\al n|^{-{1\over n}}=1~.$$
The parameters $x\delta$ and $\chi$ are in the range of analytic
functions which are defined on a Riemann surface covering the resolvent
set of the operator $h(1)$.
\medskip}
\vglue.5in

\baselineskip=.24truein

\noindent{\bf Introduction}
\medskip
The spectral properties of almost periodic operators have been
investigated extensively over the past 25 years, both in mathematics as
well as in physics. A particularly intriguing problem in this area this
author has been occupied with for some time is the question under what
conditions such operators have point spectrum and how the prevalence of
point spectrum affects the topological nature of their spectrum. In
certain cases the second part of this question appears to be intimately
related to a problem that deserves some consideration in its own right,
namely, simplistically put, under what conditions is every point in the
spectrum an eigenvalue?

In the sequel a family of almost periodic operators will be considered
for which this question has a satisfactory answer. The operators to be
considered are complex perturbations of bounded self-adjoint operators,
known as almost Mathieu operators or Harper's operators. The approach
chosen is $C^*$-algebraic. The key to the proofs are certain\break
automorphisms $\rhob$ of the irrational rotation $C^*$-algebra $\A_\al$
associated with an irrational number $\al$. If $u$ and $v$ are unitary
generators of $\A_\al$ satisfying the defining relation $uv=e^{2\pi\al
i}vu$, the operators of interest are of the form
$$h(\delta)=u+u^*+\beta(\delta v+\delta^{-1}v^*)$$
where $\beta>1$ is a fixed constant and $|\delta|>1$. The said
automorphism $\rhob$ flips a defining parameter when applied to a
slightly enlarged family of operators. This can then be used to generate
exponentially decaying eigenvectors for the operators $h(\delta)$
represented on the Hilbert space $\ell^2(\Z)$, provided $\al$ satisfies a
suitable diophantine condition. From a dynamical systems point of view,
the automorphism $\rhob$ is related to a skew translation, extending the
irrational rotation underlying the dynamics of the operators in question,
in a sense to be made precise below. It is noteworthy that as $\beta$
approaches $1$, $\rhob$ approaches an automorphism of period $4$, a
so-called ``Fourier transform''. This shows, in particular, that the
extension of the dynamical systems picture, which is so vital for the
case $\beta>1$, is no longer available in the case $\beta=1$.

Introducing the automorphism $\rhob$ and presenting a brief discussion of
the extended dynamical systems picture will be taken up in the first
paragraph. In the second paragraph the existence of exponentially
decaying eigenvectors for points in the spectrum of the operators
$h(\delta)$ will be proved. It will then be shown that a far more
specific formulation of the eigenvalue problem for the operators
$h(\delta)$ can be obtained through parametrization in a suitable Riemann
surface covering the resolvent set $\R$ of the operator $h(1)$. More
specifically, since the spectrum of $h(1)$ turns out to be a regular
compactum in the sense of potential theory, and since the spectra of the
operators $h(\delta)$ are exactly the level curves of the corresponding
conductor potential, there exists a Riemann surface $\widetilde{\R}$
covering $\R$ and an analytic function $G$ which maps $\widetilde{\R}$
onto the complement of the closed unit disk, such that
$$h(G(z))\xi={p}(z)\xi$$
has an exponentially decaying solution $\xi$ for every $z$ in $\tr$. Here
$p$ denotes the canonical mapping from $\tr$ onto $\R$. As $z$ ranges
over $\tr$, $\xi$ ranges over all possible eigenvectors for the operators
$h(\delta)$. Due to the basic $K$-theory for the $C^*$-algebra $\A_\al$,
one can see that the group of covering transformations of $\tr$ over $\R$
is infinite cyclic. Moreover, translation by one of the two generators of
this group, $\omega$ say, corresponds to shifting the eigenvector $\xi$.

The automorphism $\rhob$ gives rise to an eigenvalue problem in its own
right which is intimately interconnected with the one expounded above.
The eigenvalues are given through an analytic function $\Gamma$ on $\tr$
which has the property
$$\Gamma(\omega(z))=G(z)^2\Gamma(z)~.$$
In paragraph 3 it will be shown that the two eigenvalue problems are
essentially equivalent. To this end, the latter eigenvalue problem will
be transformed into a question regarding the kernel of certain Fredholm
operators of index zero. The problem then boils down to the question
whether these kernels are one-dimensional. With the aid of analytic
perturbation theory, it will be shown that this is indeed the case.

Finally, in paragraph 4, the case $|\delta|=1$ will be discussed. Since
$h=h(1)$ is a fixed-point of $\rhob$, no exponentially decaying
eigenvectors can be generated along the lines this was possible for the
case $|\delta|>1$. Nevertheless, the extended dynamical systems picture
shows that the eigenvalue problem for $h$ is intimately related to
similar questions about Schr\"{o}dinger type difference operators with
unbounded potentials, such as
$$n\mapsto\tan\pi(\al n^2+2\theta n+\nu)~.$$
Even though these operators are distinctly non almost periodic, they are
related to, and in fact extensions of, a family of operators which was a
focus of research in the early 1980's (let $\al=0$ and let $\theta$ be
irrational). In the restricted case ($\al=0$) it can be shown that the
said operators have pure point spectrum. This is being accomplished by
relating these operators to certain bounded operators which can be
diagonalized by solving a specific cocycle equation, and then by
translating this information back to the original context. In the
extended case ($\al\not=0$) there still exist those related bounded
operators, which are actually derived from the automorphism $\rhob$, and
which are the point of departure, rather than an auxiliary device, as it
happens to be the case when $\al=0$. But it is not possible anymore to
diagonalize these related bounded operators, due to the complications
brought about by switching from an irrational rotation to a skew
translation extending it.

\bigskip
\bigskip
\item{\bf 1.}{\bf The automorphism $\rho_\beta$}
\medskip
\noindent For an irrational number $\al$ we consider the $C^*$-algebra $\A=\A_\al$
generated by two unitaries $u$ and $v$ satisfying the relation
$uv=\lam^2vu$, where $\lam=e^{\pi\al i}$. For $\be\not=1$, let
$$\eqalign{
\rho_\beta(u)&=vuv(uv+\beta)^{-1}(v^*u^*+\beta)\cr
\rho_\beta(v)&=v(uv+\beta)^{-1}(v^*u^*+\beta)~.\cr}$$
Then $\rhob(u)$ and $\rhob(v)$ are unitaries satisfying again the relation
$$\rhob(u)\rhob(v)=\lam^2\rhob(v)\rhob(u)~.$$
Therefore, $\rhob$ extends to an automorphism of $\A$ which we will also
denote by $\rhob$. The significance of this automorphism for what is to
follow rests with the identities
$$\rhob(u+\be v)=u^*+\be v~~,~~\rhob(u^*+\be v^*)=u+\be v^*~.\leqno(1.1)$$
Let $GL(2,\Z)$ be the group of $2\times 2$ matrices with integer entries
and a determinant of modulus 1. The assignment
$$
\left.
{\eqalign{w_{m,n}&\to\hbox{$w_{pq}$}\cr
\left({p\atop q}\right)&=\hbox{$A\left({m\atop n}\right)$}\cr}}
\right\}\quad A\in GL(2,\Z)~,
$$
where $w_{pq}=\lam^{-pq}u^pv^q$, is known to define a linear isometry of
$\A$. This isometry is an automorphism if $\det A=1$, and it is an
antiautomorphism if $\det A=-1$. The mapping assigning to every $A\in
GL(2,\Z)$ the corresponding isometry is a faithful homomorphism from the
group $GL(2,\Z)$ into the group of isometries of $\A$. For convenience we
will denote matrices in $GL(2,\Z)$ and their corresponding isometries by
the same symbol. For instance $\left({0~1\atop{1~0}}\right)$ represents
the antiautomorphism $u\mapsto v$, $v\mapsto u$.

Returning to the automorphism $\rhob$ we are now going to list a number of
useful identities.

$$
\left.
{\eqalign{
{}&\rho_\beta=\left({1~0\atop{1~1}}\right)\circ\rho^{(0)}_{\beta}\circ
\left({1~0\atop{1~1}}\right)~,\hbox{where}\cr
{}&\rho^{(0)}_{\beta}(u)=u\cr {}&\rho^{(0)}_{\beta}(v)=v(\lambda
u+\beta)^{-1}(\ovl{\lambda}u^*+\beta)\cr}}\right\}\leqno(1.2)
$$

$$\rho_{\be^{-1}}=\pmatrix{0&-1\cr 1&~~0}\circ\rho^{-1}_\beta\circ
\pmatrix{0&-1\cr 1&~~0}\leqno(1.3)$$
This identity says essentially that $\rho_{\be^{-1}}$ and
$\rho^{-1}_\beta\circ\({-1~~~0\atop ~~0~-1}\)$ are conjugates of each
other via a ``Fourier transform'' (also known as ``duality'').
$$\pmatrix{-1&~~0\cr ~~0&-1}\circ\rho_\beta=\rho_\beta\circ\pmatrix{-1&~~0\cr
~~0&-1}\leqno(1.4)$$ If $\beta$ approaches $1$ then $\rho^{(0)}_\beta$
approaches $\({1~-1\atop 0~~~1}\)$ on $w_{pq}$, hence
$\rho_\be$\hfill\break approaches
$$\pmatrix{1&0\cr 1&1}\circ\pmatrix{1&-1\cr 0&~~1}\circ
\pmatrix{1&0\cr 1&1}=\pmatrix{0&-1\cr 1&~~0}$$
on $w_{pq}$. Thus $\rhob$ approaches a ``Fourier transform''. Throughout
in the subsequent discussion we will limit our attention to irrational
numbers $\al$ only which satisfy a diophantine condition.
$$\lim\limits_{n\to\infty}|\sin\pi\al n|^{-1}=1\leqno(1.5)$$
For such numbers $\al$ and $\be>1$ the automorphism $\rho^{(0)}_\beta$
becomes an inner automorphism. Namely,
$$
\left\{
{\eqalign{{}&\rho^{(0)}_\beta(a)=e^{ig(u)}ae^{-ig(u)}~, a\in{\A}~;
\hbox{or}~\rho^{(0)}_{\beta}= Ad(e^{ig(u)})~,\cr
{}&\hbox{where}\cr
{}&g(z)=\sum\limits^{\infty}_{n=1}{(-1)^n\over{n}}(\sin\pi\al
n)^{-1}\beta^{-n}(z^n+z^{-n})~.\cr}}
\right.\leqno(1.6)
$$

We turn now to a dynamical systems interpretation of the automorphisms
$\rhob$ in the framework of $C^*$-algebras, under the assumption that
(1.5) holds. Consider the crossed product $\B$ of $\A$ by the
automorphism $\({1~0\atop 1~1}\)$
$$\B=\A\otimes_{{\left(\hbox{$1~0\atop{1~1}$}\right)}}\Z~.$$
This $C^*$-algebra is generated by three unitaries $u$, $v$ and $w$
satisfying the defining relations
$$
\left.
{\eqalign{w^*uw&=\ovl{\lambda}uv\cr vw&=wv\cr
uv&=\lambda^2vu\cr}}\right\}
$$
Since $\rhob=\({1~0\atop 1~1}\)\circ Ad(e^{ig(u)})\circ\({1~0\atop
1~1}\)$, $\rhob$ extends to an inner automorphism of $\B$. So $\B$
provides a natural framework where all our manipulations so far take
place. Manipulating the first of those relations we get
$$uwu^*=\lambda vw~.$$
This, in conjunction with the other two relations, suggests that $\B$ can
also be realized as the crossed-product of the $C^*$-algebra of
continuous functions on the two dimensional torus $C(\T^2)$ by a skew
translation followed by a translation. More specifically, let
$$\phi(x,y)=\lambda(\lambda x,xy)~,\quad x,y\in\T^2~.$$
In dynamical systems theory it is a well known fact that $\phi$ is a
uniquely ergodic homeomorphism of $\T^2$, the unique invariant
probability measure being the Haar measure on the compact group $\T^2$,
which has of course full support. It is then a well known fact in
$C^*$-algebra theory, that the crossed-product
$$\tb=C(\T^2)\otimes_{\phi}\Z$$
is a simple $C^*$-algebra (i.e., it has no non-trivial ideals) with a
unique tracial state\break extending the Haar measure on $C(\T^2)$.
Letting $\wt u$ be the unitary in $\tb$ corresponding to $\phi$, $\wt v$
the projection from $\T^2$ onto the first component, and $\wt w$ the
projection from $\T^2$ onto the second component, then it is easily seen
that $\wt u$, $\wt v$ and $\wt w$ satisfy the same three relations stated
above for $u$, $v$ and $w$. Using this information it is not difficult to
see that the assignments ${\wt u}\mapsto u$, ${\wt v}\mapsto v$, ${\wt
w}\mapsto w$ extend to an isomorphism from $\tb$ onto $\B$.

\vfill\eject

\item{\bf 2.}{\bf Spectrum and point spectrum for a family of
non self-adjoint almost\hfill\break periodic operators}
\medskip

\noindent In this paragraph we will employ the automorphism $\rhob$ to
investigate the spectrum of the operators
$$h(\delta)=u+u^*+\beta(\delta v+\delta^{-1}v^*)$$
for $|\delta|>1$, provided that $\be>1$ and $\al$ satisfies the property
(1.5). To this end we consider the extended family
$$h_\gam(\delta)=\gam u+\gam^{-1}u^*+\be(\de v+\de^{-1}v^*)~,$$
where $\beta^{-1}\delta^{-1}<|\gam|<\be\de$.

The assignments
$$\left\{\eqalign{(u\xi)_n&=\xi_{n-1}\cr
(v\xi)_n&=\ovl{\lam}^{2n}\xi_n\cr
(w\xi)_n&=\lam^{n^2}\xi_n}\right.\leqno(2.1)$$ define linear operators on
the vector space $\C^\infty$ of (two-sided) sequences $\xi$. When
restricted to square summable sequences, these assignments extend to a
(faithful) representation of the $C^*$-algebra $\B$ on the Hilbert space
$\ell^2(\Z)$ introduced in paragraph 1. For $\be^{-1}<|\gam|<\be$, let
$$k_\gam=we^{ig(\gam u)}w~,\quad k=k_1~.$$
Then (1.1), (1.2) and (1.6) yield
$$h_{\gam\del^{-1}}(\del)k_\gam=k_\gam h_{\gam\del}(\del)~.\leqno(2.2)$$

For any complex number $x$, let
$$(D_x\xi)_n=x^n\xi~,\quad \xi\in\C^\infty~.$$

\noindent{\bf 2.1 Lemma.} Let $|\del|>1$; $\be^{-1}<|\gam|<1$ or
$1<|\gam|<\be$.

\noindent{$(+)$} Suppose that
$$h_{\gam\del}(\del)\eta=z\eta~\hbox{for some}~z\in\C~, \eta\in\ell^\infty(\Z)~.$$
Then
$$h(\del)\xi=z\xi~,\quad \limsup_{n\to-\infty}|{\xi_n}|^{1/|n|}
\le|\gam|\del^{-1}~,~\limsup_{n\to\infty}|\xi_n|^{1/n}\le|\gam|^{-1}\del^{-1}$$
where $\xi=D_{\gam^{-1}\del^{-1}}\eta$.

\vfill\eject
\noindent{$(-)$} Suppose that
$$h_{\gam\del^{-1}}(\del)\eta=z\eta~~~\hbox{for some}~z\in\C,
~~\eta\in\ell^\infty(\Z)~.$$
Then
$$h(\del)\xi=z\xi~,\limsup_{n\to-\infty}|\xi_n|^{1/|n|}\le|\gam|^{-1}\del^{-1}~,~
\limsup_{n\to\infty}|\xi_n|^{1/n}\le|\gam|\del^{-1}~,$$
where $\xi=D_{\gam^{-1}\del}\eta$.

\bigskip
\noindent{\bf Proof:} We will deal with the case $(+)$ only. The case $(-)$
can be handled in a similar fashion.

The identity (2.2) yields
$$h_{\gam\del^{-1}}k_\gam\eta=zk_\gam\eta~,~~~\hbox{where}~~~
k_\gam\eta\in\ell^\infty(\Z)~.$$ Let $\txi=D_{\gam^{-1}\del}\eta$. Then
$$h(\del)\txi=z\txi~,~~\limsup_{n\to-\infty}|\txi_n|^{1/|n|}\le|\gam|\del^{-1}~,~
\limsup_{n\to\infty}|\txi_n|^{1/n}\le|\gam|^{-1}\del~.$$
On the other hand,
$$h(\del)\xi=z\xi~,~~\limsup_{n\to-\infty}|{\wt \xi}_n|^{1/|n|}\le|\gam|\del~,~
\limsup_{n\to\infty}|\xi_n|^{1/n}\le|\gam|^{-1}\del^{-1}~.$$
Therefore

$$\limsup_{n\to-\infty}\left|\xi_n\txi_{n+1}\right|^{1/|n|}\le|\gam|^2~,~
\limsup_{n\to\infty}\left|\xi_n\txi_{n+1}\right|\le|\gam|^{-2}$$
and
$$\limsup_{n\to-\infty}\left|\xi_{n+1}\txi_n\right|^{1/|n|}\le|\gam|^2~,~
\limsup_{n\to\infty}\left|\xi_{n+1}\txi_n\right|\le|\gam|^{-2}~.$$
Since $|\gam|\not=1$ by assumption, this entails
$$\limsup_{n\to-\infty}\left(\xi_n\txi_{n+1}-\xi_{n+1}\txi_n\right)=0$$
or
$$~~\limsup_{n\to\infty}\left(\xi_n\txi_{n+1}-\xi_{n+1}\txi_n\right)=0~.$$
Either way, it follows that $\xi$ and $\txi$ are linearly dependent. For,
if this were not the case, the expression following the $\limsup$ in the
last two identities would have to be constant and non-zero for all
$n\in\Z$. Therefore, observing that $\del>1$ by assumption, we conclude
$$\lim_{n\to-\infty}|\xi_n|^{1/|n|}\le|\gam|\delta^{-1}~,
~~\lim_{n\to\infty}|\xi_n|^{1/n}\le|\gam|^{-1}\del^{-1}~,$$ as
claimed.\qed

\noindent{\bf 2.2 Lemma.} Let $\I$ be a set of real numbers and let $b>0$,
$d>0$. Suppose that $\I$ has a non-empty intersection with at least one
of the two open intervals $(d-b, d)$ or $(d, d+b)$. Suppose in addition
that $\I$ has the following properties:
\item{$(+)$} $0<|t|<b$ and $t+d\in\I$ implies $(t-d, t+d)\subset\I$.
\item{$(-)$} $0<|t|<b$ and $t-d\in\I$ implies $(t-d, t+d)\subset\I$.

\noindent Then $(-d-b, d+b)\subset\I$.

\bigskip
\noindent{\bf Proof:} Consider the case that
$\I\cap(d-b,d)\not=\emptyset$. Then there exists a $c_0\in(-b,0)$ such
that $c_0+d\in\I$.

Suppose first that $b\le d$.

Let $s>-d-b$, but close to $-d-b$. Since $(+)$ ensures that $(c_0-d,
c_0+d)\subset\I$, we can find $c_1\in(-b,0)$ such that $c_1+d\in\I$ and
$c_1-d<s$. Again, $(+)$ ensures that $(c_1-d, c_1+d)\subset\I$. Now let
$t<d+b$ but close to $d+b$. Then we can find $c_2\in(0,b)$ such that
$c_2-d\in(c_1-d,c_1+d)\subset\I$ and $c_1+d>t$. Then $(-)$ ensures that
$(c_2-d, c_2+d)\subset\I$. By construction $(s,t)\subset(c_1-d,
c_1+d)\cup(c_2-d, c_2+d)\subset\I$. Since $s$ and $t$ can be chosen
arbitrarily close to $-d-b$ and $d+b$, respectively, we conclude that
$(-d-b, d+b)\subset\I$.

Now suppose that $b>d$.

Let $s>-d-b$, but close to $-d-b$. Through induction we can generate a
(possibly empty) chain $c_1>\ldots>c_n$ such that $c_0>c_1$ if $n>0$;
$c_k\in(-b,0)$,\hfill\break $c_k\in(c_{k-1}-d, c_{k-1}+d)$, for $1\le
k\le n$, and $c_n-d<-b$. Repeated applications of $(+)$ show that
$(c_k-d, c_k+d)\subset\I$ for $k=0,\ldots,n$. Now choose
$c_{n+1}\in(-b,c_n)$ such that $c_{n+1}-d<s$. Once again, $(+)$ ensures
that $(c_{n+1}-d, c_{n+1}+d)\subset\I$. By construction we have
$s\in(c_{n+1}-d,c_0+d)\subset\I$.

Next, let $t<d+b$ but close to $d+b$. Using $(-)$ instead of $(+)$, we
can construct in the same fashion a chain $c_0<\tc_1<\ldots<\tc_{m+1}$
such that $t\in(c_0-d, c_{m+1}+d)\subset\I$. This means $(s,t)\subset\I$,
and again we conclude that $(-d-b, d+b)\subset\I$.

The case that $\I\cap(d,d+b)\not=\emptyset$ can be handled in a similar
fashion.\qed

\noindent{\bf 2.3 Lemma.} Let $|\delta|>1$; $\beta^{-1}<|\gam_0|<1$ or
$1<|\gam_0|<\beta$. Suppose that
$$h_{\gam_0\delta^{\pm 1}}(\delta)\eta=z\eta~~\hbox{for some}~~z\in\C~,
~~\eta\in\ell^{\infty}(\Z)~.$$ Then
$$h(\delta)\xi=z\xi~~,~~\limsup_{|n|\to\infty}
|\xi_n|^{1/|n|}\le\beta^{-1}|\delta|^{-1}~,$$ where
$\xi=D_{\gam_0^{-1}\delta\mp 1}\eta$.
\bigskip
\noindent{\bf Proof:} Let
$$\displaylines{b=\log\beta~,~d=\log|\delta|~,~c_0=\log|\gam_0|~,\cr
\I=\left\{\log|\gam|~\big|~D_\gam\xi\in\ell^\infty(\Z)\right\}~.\cr}$$
Then $(+)$ and $(-)$ in Lemma 2.1 translate into the namesake properties
of Lemma 2.2, which then yields the desired conclusion.\qed

\bigskip
Let $Sp_0(h_\gam(\del))$ be the spectrum of $h_\gam(\del)$ considered as
a bounded linear operator on the Banach space $c_0(\Z)$ of bounded
two-sided sequences which vanish at infinity.

\bigskip
\noindent{\bf 2.4 Lemma.} For every $z\in Sp_0(h_\gam(\del))$ there
exists $x\in\T$ and $\eta\in\ell^\infty(\Z)\backslash\{0\}$ such that
$$h_\gam(x\del)\xi=z\xi~.$$

\noindent{\bf Proof:} We choose an approximate eigenvector for $z$ in
$c_0(\Z)$, $\eta^{(1)},\eta^{(2)},\ldots$ in $c_0(\Z)$,
$\|\eta^{(m)}\|_\infty=1$,
$$\lim_{m\to\infty}\|(h_\gam(\del)-z))\eta^{(m)}\|_\infty=0~.$$
For each $m$ there is a $j_m$ such that
$\left|\eta^{(m)}_{j_m}\right|\ge{1\over 2}$. Let
$\xi^{(m)}=u^{-j_m}\eta^{(m)}$. Then
$$\lim_{m\to\infty}\|(h_\gam(x_m\del)-z)\xi^{(m)}\|_\infty=0~,$$
where $x_m=\ovl{\lambda}^{2j_m}$. Also,
$$\left|\xi_0^{(m)}\right|\ge {1\over 2}~.$$
Since the unit ball of $\ell^\infty(\Z)$ is weakly compact and
metrizable, there exists a subsequence of $\xi^{(1)},\xi^{(2)},\ldots$
which converges weakly to some $\xi\in\ell^\infty(\Z)$. Moreover, we can
arrange for the corresponding subsequence of $x_1,x_2,\ldots$ to converge
to some $x\in\sprod$. Since $|\xi_0|\ge{1\over 2}$, $\xi$ is non-zero. By
construction
$$h_\gam(x\del)\xi=z\xi\qeddis

\noindent{\bf 2.5 Lemma.} $Sp_0(h_\gam(\del))=Sp_0(h_{|\gam|}(|\del|))$.

\bigskip
\Proof:  Since
$h_\gam(\del)=D_{\gam/|\gam|}h_{|\gam|}(\del)D^{-1}_{\gam/|\gam|}$, we
have
$$Sp_0\(h_\gam(\del)\)=Sp_0\(h_{|\gam|}(\del)\)~.$$
Now let $z\in Sp_0(h_\gam(\del))$ and $x\in\T$. Then there exists an
approximate eigenvector for $z$ in $c_0(\Z)$:
$$\eta^{(1)},\eta^{(2)},\ldots ;~~~ \|\eta^{(m)}\|_\infty=1~,$$
$$\lim_{m\to\infty}\|(h_\gam(\del)-z)\eta^{(m)}\|_{\infty}=0~.$$
For each $m$ choose $j_m\in\Z$ such that $\lam^{2j_m}\to x$. Then
$u^{j_1}\eta^{(1)},u^{j_2}\eta^{(2)},\ldots$ is seen to be an approximate
eigenvector of $h_\gam(x\del)$ for $z$, which entails $z\in
Sp_0(h_\gam(x\del))$.\qed

\bigskip
\noindent{\bf 2.6 Lemma.} Let $|\del|>1$,
$|\del|\beta^{-1}<|\gam|<|\del|\beta$. Then
$$Sp_0\(h_\gam(\del)\)=Sp_0\(h_{\gam^{-1}}(\del)\)~.$$

\bigskip
\Proof: Let $z\in Sp_0(h_\gam(\del))$. By Lemma 2.4 there exist
$x\in\sprod$, $\xi\in\ell^\infty(\Z)\backslash\{0\}$ such that
$$h_\gam(x\del)\xi=z\xi~.$$
Lemma 2.3 implies that $\xi$ and $\eta=D_{\gam^{-2}}\xi$ are in
$c_0(\Z)$. But
$$h_{\gam^{-1}}(x\del)\eta=z\eta~.$$
Hence
$$z\in Sp_0(h_{\gam^{-1}}(x\del))~,$$
which, by Lemma 2.5, entails
$$z\in Sp_0(h_{\gam^{-1}}(\del))~.$$
so we have
$$Sp_0(h_\gam(\del))\subset Sp_0(h_{\gam^{-1}}(\del))~.$$
The opposite inclusion is shown in exactly the same way.\qed

\noindent{\bf 2.7 Lemma.} Let $|\del|>1$,
$|\del|\beta^{-1}<|\gam|<|\del|\beta$. Then
$$Sp(h_\gam(\del))=Sp_0(h_\gam(\del))~,$$
where $Sp(h_\gam(\del))$ denotes the spectrum of $h_\gam(\del)$
considered as an operator on the Hilbert space $\ell^2(\Z)$.

\Proof: First we are going to show that $Sp_0(h_\gam(\del))\subset
Sp(h_\gam(\del))$. Let $z\in\C\backslash Sp(h_\gam(\del))$, and let
$$\left(h_\gam(\del)-z\right)^{-1}=\sum^{\infty}_{p,q=-\infty}c_{pq}(z)~w_{pq}$$
be the ``Fourier expansion'' of $(h_\gam(\del)-z)^{-1}$. Since this
expansion decays exponentially in $p$ and $q$, it defines a bounded
linear operator on $c_0(\Z)$, which is an inverse of $h_\gam(\del)-z$ on
$c_0(\Z)$. Hence $z\in\C\backslash Sp_0(h(\del))$. So,
$Sp_0(h_\gam(\del))\subset Sp(h_\gam(\del))$.

To establish the opposite inclusion, let $z\in\C\backslash
Sp_0(h_\gam(\del))$. Then $h_\gam(\del)-z$ has an inverse $S$ on
$c_0(\Z)$. Since $\spo(h_\gam(\del))=\spoh$ by Lemma 2.6, $\hgodl-z$ also
has an inverse on $c_0(\Z)$, which we denote by $T$. Let $(\hgodl-z)^t$
and $T^t$ be the transposed operators of $\hgodl-z$ and $T$,
respectively, on the dual Banach space of $c_0(\Z)$, which happens to be
$\ell'(\Z)$. Since $(\hgodl-z)^t$ equals the restriction of $\hgd-z$ on
$\ell^1(\Z)\subset c_0(\Z)$ and $T^t$ is an inverse of $(\hgodl-z)^t$ on
$\ell^1(\Z)$, $T^t$ must be the restriction of $S$ on $\ell^1(\Z)$. Now
consider the double transposition $(\hgd-z)^{tt}$ and $S^{tt}$ of
$\hgd-z$ and $S$, respectively, on the dual Banach space of $\ell^1(\Z)$,
which happens to be $\ell^\infty(\Z)$. Then $(\hgd-z)^{tt}$ restricted on
$c_0(\Z)$ equals $\hgd-z$, while $S^{tt}$ restricted on $c_0(\Z)$ equals
$S$. Thus, dispensing with the double $t$'s, we can summarize the
situation as follows:

Considering the linearly embedded Banach spaces
$\ell^1(\Z)\subset\ell^2(\Z)\subset\ell^\infty(\Z)$, we have two linear
operators on $\ell^\infty(\Z)$, namely, $\hgd-z$ and $S$, which are
inverses of each other. Both map $\ell^1(\Z)$ into $\ell^1(\Z)$, and both
are continuous on $\ell^\infty(\Z)$ as well as on $\ell^1(\Z)$ with
regard to the respective Banach space norms. In addition, $\hgd-z$ maps
$\ell^2(\Z)$ (continuously) into $\ell^2(\Z)$. On account of the
Riesz-Thorin theorem ([Kn], Ch. IV), we conclude that $S$ maps
$\ell^2(\Z)$ into $\ell^2(\Z)$ and that $S$ restricted on $\ell^2(\Z)$ is
continuous with respect to the Hilbert space norm. Moreover, $S$ and
$\hgd-z$, restricted on $\ell^2(\Z)$, are inverses of each other. Hence
$z\in\C\backslash Sp(h(\del))$, and therefore $Sp(h_\gam(\del))\subset
Sp_0(h_\gam(\del))$.\qed

\noindent{\bf 2.8 Lemma.} Let $|\del|>1$. Then $\sphd$ is constant for
$|\del|^{-1}\beta^{-1}<|\gam|<|\del|\beta$.

\bigskip
\Proof: The Lemmas 2.3, 2.4 and 2.7 show that $\sphd$ is constant for
$|\del|^{-1}\beta^{-1}<|\gam|<|\del|\beta$, with the possible exception
when $|\gam|=|\del|$ and $|\gam|=|\del|^{-1}$. Let $Sp(\del)$ denote the
constant spectrum covered by those cases. The lemmas listed also show
that
$$Sp(\del)\subset \sphd~~\hbox{for}
~~|\del|^{-1}\beta^{-1}<|\gam|<|\del|\beta~.$$ Now consider the ``Fourier
expansion'' of the resolvent of $h(\del)$
$$(h(\del)-z)^{-1}=\sum^\infty_{p,q=-\infty}c_{pq}(z)w_{pq}~.$$
We have $c_{pq}(z)=c_{|p|,q}(z)$. From [R2], paragraph 4, we know that
$$(\hgd-z)^{-1}=\sum^{\infty}_{p,q=-\infty}
\gam^pc_{pq}(z)w_{pq}~~\hbox{for}~~z\in\C\backslash Sp(\hgd)~.$$

These two facts taken together show that $Sp(\hgd)$ increases as
$\bigl|\log|\gam|\bigr|$ increases. In conjunction with the inclusion
established above, this shows that $Sp(\hgd)=Sp(\del)$ for
$|\del|^{-1}\beta^{-1}<|\gam|<|\del|\beta$, as claimed.\qed

\noindent{\bf Remark.} The proof of Lemma 2.8 looks a bit twisted. Since
it is not this author's aspiration to state and prove the relevant facts
in their utmost generality, he feels that the line of reasoning chosen in
the context of the objective in this exposition is somewhat adequate.
However, it should be mentioned that the claim of Lemma 2.8 is valid for
arbitrary irrational numbers $\al$. Here is the brief sketch of a proof.

Consider the set $\A^{(\omega)}$ of analytic elements in $\A$. An element
is called analytic if it has an exponentially decaying ``Fourier
expansion''. $\A^{(\omega)}$ is seen to be a $*$-subalgebra of $\A$. Even
though the concept of the operator $k_\gam$, as defined at the beginning
of this paragraph, cannot be extended to the general setup, the concept
of the algebra automorphism $Ad(k_\gam)$ can. To this end, one needs to
show that
$$Sp\((e^{\pi\al i}\gam u+\beta)^{-1}(e^{-\pi\al i}\gam^{-1}u+\beta)\)
=\T,~~\hbox{whenever}~~\beta^{-1}<|\gam|<\beta~,$$
for all irrational numbers $\al$. This can be fairly easily proved by
exploiting the spectral radius formula as well as the ergodicity of the
irrational rotation on $\T$ associated with $\al$, in conjunction with
Ascoli's theorem. One can use this information to show that the
assignments
$$\eqalign{u&\longmapsto u\cr v&\longmapsto
v(e^{\pi\al i}\gam+\beta)^{-1} (e^{-\pi\al i}\gam^{-1}u^*+\beta)}$$ yield
an algebra automorphism $\rho_{(\gam)}^{(0)}$ of $\A^{(\omega)}$ (which
does of course not preserve the involution $*$ unless $|\gam|=1$).

Multiplying $\rho_{(\gam)}^{(0)}$ from both sides by the automorphism
$\({1~0\atop 1~1}\)$ yields an automorphism $\rho_{(\gam)}$ of
$\A^{(\omega)}$ which has the property
$$\rho_{(\gam)}\(h_{\gam\del}(\del)\)=h_{\gam\del^{-1}}(\del)~~
\hbox{for}~~\beta^{-1}<|\gam|<\beta~.$$
This identity extends (2.2). Since an element in $\A^{(\omega)}$ is
invertible in $\A^{(\omega)}$ if and only if it is invertible in $\A$,
$\rho_{(\gam)}$ preserves the spectrum of an element in $\A^{(\omega)}$.

Now consider the resolvent $(h_\gam(\del)-z)^{-1}$ of $\hgd$. There is a
power series
$$\sum^{\infty}_{p,q=0}c_{pq}(z)x^py^q$$
in two variables such that
$$\(\hgd-z\)^{-1}=\sum^{\infty}_{p,q=-\infty}\tc_{pq}(z)\gam^p\del^q w_{pq}~,$$
where throughout in each of the four $(p,q)$-quadrants either
$$\tc_{pq}(z)=c_{|p|,|q|}(z)~~\hbox{ or }~~ \tc_{pq}(z)=0$$
(cf. [R2], paragraph 4). This power series depends on $z$ only and the
type of the component in the resolvent set of $\hgd$ to which $z$ belongs
(i.e. bounded component vs. unbounded component). Since $\rho_{(\gam)}$
preserves the spectrum of any element in $\A^{(\omega)}$, one can show
that the domain of convergence of the said power series, which is known
to be a logarithmically convex complete Reinhardt domain, must contain
the polydisc
$$\{x\in\C~\big|~|x|\le d\cdot\beta\}\times\{y\in\C~\big|~|y|\le d\}~,$$
where $d=\max\{|\del|,|\del|^{-1}\}$. This is valid whenever $|\del|>0$
and $d^{-1}\beta^{-1}<|\gam|<d\beta$. It follows immediately that
$Sp(\hgd)$ is constant for $d^{-1}\beta^{-1}<|\gam|<d\beta$. (Notice that
no exception needs to be made for the case $|\del|=1$.)\qed

\bigskip
We are now in a position to record a preliminary conclusion to our
discussion regarding the existence of eigenvectors.

\bigskip
\noindent{\bf 2.9 Proposition.} For every $\del\in\C\backslash(\T\cup\{0\})$
and every $z\in Sp(h(\del))$ there exist $x\in\T$,
$\xi\in\ell^2(\Z)\bsl\{0\}$ such that
$$h(x\del)\xi=z\xi~,\quad \lim_{|n|\to\infty}\left|\xi_n\right|^{1/|n|}
\le\min\left\{|\del|\beta^{-1},|\del|^{-1}\beta^{-1}\right\}~.$$

\Proof: If $|\del|>1$, the claim follows from Lemmas 2.3, 2.4, 2.7 and
2.8. Since
$${\J}h(\del){\J}=h(\del^{-1})~,~~\hbox{where}~({\J}\xi)_n=\xi_{-n}~,$$
the case $0<|\del|<1$ can be reduced to the case $|\del|>1$.\qed

\bigskip
The translation by $x$ in Proposition 2.9, and hence $\xi$, is
(essentially) uniquely determined by $z$.

\bigskip
\noindent{\bf 2.10 Proposition.} If $\del>1$; $x,y\in\T$ and
$$h(x\del)\xi=z\xi~, h(y\del)\eta=z\eta~; \xi,\eta\in\ell^2(\Z)~,$$
then $y=\lam^{2\ell}x$ for some $\ell\in\Z$ and $u^\ell\xi$, $\eta$ are
linearly dependent.

\bigskip
\Proof:  We have
$$h_x(x\del)D_x\xi=zD_x\xi~,
\quad h_{\ovl y}(y\del)D_{\ovl y}\eta=zD_{\ovl y}\eta~.$$
So, if we let
$$\varphi(a)=\langle a D_x\xi, D_y{\ovl\eta}\rangle~,\quad a\in\A~,$$
then
$$\varphi(ah_x(x\del))=\varphi(h_y(y\del)a)=z\varphi(a)~.$$
Now let
$$\varphi_{pq}=s^{p+q}\del^{-q}\varphi(w_{pq})~,
\quad\hbox{where}~~s^2={\ovl x}{\ovl y}~.$$
Then the double sequence $\left\{\varphi_{pq}\right\}$ is seen to solve
the system of difference equations
$${\eqalign{\cos(\pi\al q+\theta)(X_{p-1,q}+X_{p+1,q})
&+\beta\cos(\pi\al p+\theta)(X_{p,q-1}+X_{p,q+1})=zX_{pq}\cr
\sin(\pi\al q+\theta)(X_{p-1,q}-X_{p+1,q})
&-\beta\sin(\pi\al p+\theta)(X_{p,q-1}-X_{p,q+1})=0~,}}\leqno(*)$$ where
$e^{i\theta}=s$. Starting over again with $\ovl\xi$ and $\ovl\eta$ in
place of $\xi$ and $\eta$, respectively, we get
$$h_x(x\del^{-1})D_x\oxi=\z D_x\oxi~,~~h_{\oy}
(y\del^{-1})D_{\oy}\oeta=\z D_{\oy}\oeta~.$$ Let
$$\psi(a)=\langle a D_x\oxi, D_y\eta\rangle~,\quad a\in\A~.$$
Now let
$$\psi_{pq}=s^{p+q}\delta^q\psi(w_{pq})~,\quad\hbox{$s$ as before.}$$
Replacing $z$ by $\z$ in $(*)$ it is seen that the double sequence
$\{\psi_{pq}\}$ solves $(*)$. So both, $\{\varphi_{pq}\}$ and
$\{\ovl{\psi}_{pq}\}$ solve $(*)$ for the same parameter $z$. Moreover,
the sequence
$$\left\{\varphi_{pp}\psi_{p+1,p+1}-\varphi_{p+1,p+1}\psi_{pp}\right\}$$
is bounded.

Now suppose that the claim of the proposition is not true. Then it
follows from appendix A1 that $\{\varphi_{pq}\}$ and
$\{\ovl{\psi}_{pq}\}$ are linearly dependent. This however implies that
$\{\varphi_{pq}\}$ decays exponentially uniformly in $p$ of order at
least $\del^{-2}$ as $q\to\infty$. Since the sequence
$\{\psi_{pq}\}_{q\in{\bs Z}}$ is almost periodic for every $p\in\Z$, and
hence does not approach zero as $q\to\infty$, we have reached a
contradiction.\qed

\bigskip
Our next objective is to show that the spectrum of $h=u^*+u+\beta(v+v^*)$
is a regular compactum in the sense of potential theory. In view of [R3],
Theorem 2.2 we are going to prove a stronger statement. In preparation of
this, we need the following.

\bigskip
\noindent{\bf 2.11 Lemma.} $Sp(h)\subset\C\bsl Sp(h(\del))$ whenever
$|\del|\not=1$.

\bigskip
\Proof: Suppose this were not true. Then there exists $\chi\in Sp(h)$
and $\del_0>1$ such that $\chi\in Sp(h(\del))$ for
$\del^{-1}_0\le\del\le\del_0$. According to Proposition 2.9, for any
$\del\in[\del^{-1}_0,\del_0]\bsl\{1\}$, there exist $x\in\T$ and
$\xi\in\ell^2(\Z)\bsl\{0\}$ such that
$$h(x\del)\xi=\chi\xi~,
\quad \limsup_{|n|\to\infty}|\xi|^{1/|n|}\le\del^{-1}\beta^{-1}~,$$

For every $y\in\T$, let
$$\eta^{(y)}_n=\sum^{\infty}_{p=-\infty}y^p x^{-n}\del^{-n}
\langle w_{pn}\xi,\J\xi\rangle
w_{pn}e^{(0)}~,$$ where ($\J$ as in 2.9),
$$e^{(0)}_n=\cases{1~,&$n=0$\cr 0~,&$n\not=0~.$}$$
Let $\h(y)=u+u^*+\beta^{-1}(yv+\oy v^*)$. Then
$$\h(y)\eta^{(y)}=\chi\eta^{(y)}~,\quad \limsup_{n\to\infty}
\left|\eta^{(y)}_n\right|^{1/n}\le\del^{-1}~, \quad \limsup_{n\to-\infty}
\left|\eta^{(y)}_n\right|^{1/|n|}\le\del~.$$
Exempting for every $\del\in\left[\del^{-1}_0,\del_0\right]\bsl\{1\}$ a
possibly non-empty, but countable subset of values for the parameter $y$,
for which $\eta^{(y)}$ might be zero, this implies the following: There
is a countable subset $M\subset\T$ and a dense countable subset
$N\subset\left[\del^{-1}_0,\del_0\right]\bsl\{1\}$ such that for every
$y\in\T\bsl M$ and every $\del\in N$ there exists
$\eta\in\C^\infty\bsl\{0\}$ such that
$$\h(y)\eta=\chi\eta~,
\quad\limsup_{n\to\infty}\left|\eta_n\right|^{1/n}\le\del^{-1}~,
\quad\limsup_{n\to-\infty}\left|\eta_n\right|^{1/|n|}\le\del~.$$
This in turn entails that for every $y\in\T\bsl M$, the operator $\h(y)$
has an exponentially decaying eigenvector for $\chi$. Since $\h(y)$ is
known to have no eigenvalues for any $y\in\T$, we have reached a
contradiction.\qed

\bigskip
Letting $\mu$ be the probability measure of $Sp(h)$ obtained through
restricting the canonical trace $\tau$ of $\A$ on the $C^*$-algebra
generated by $h$, Lemma 2.11 entails, on account of [R3], Theorem 2.2:

\bigskip
\noindent{\bf 2.12 Proposition.} $Sp(h(\del))=\{z\in\C\mid\int\log|z-s|d\mu(s)
=\log(\beta\del)\}~\rm{for}~\del\ge 1$. In particular $Sp(h)$ is a
regular compactum and $\mu$ is its equilibrium distribution.

\bigskip
\noindent{\bf 2.13 Corollary.} $Sp(h)$ is not connected.

\bigskip
\Proof: Consider the moments $\tau(h^n)$ of $h$. While $\tau(h^{2n-1})=0$
for every $n\in\N$, elementary calculations show that
$$\tau(h^2)=2\beta^2+2~,\quad \tau(h^4)=6\beta^4+(24+16\cos 2\pi\al)\beta^2+6~.$$
Now consider
$v+v^*=\lim\limits_{\beta\to\infty}\beta^{-1}(u+u^*+\beta(v+v^*))$. This
element has a connected spectrum and $\tau$ restricted to the
$C^*$-algebra generated by $v+v^*$ is nothing but the equilibrium
distribution for $Sp(v+v^*)=[-2,2]$. Again we have
$\tau(v+v^*)^{2n-1})=0$ for every $n\in\N$. Moreover,
$$\tau((v+v^*)^2)=2~,\quad \tau((v+v^*)^4)=6~.$$
In order to identify possible values for $\beta$ for which $Sp(h)$ is
connected, we have to solve the equations
$$2x=\tau(h^2)~,\quad 6x^2=\tau(h^4)~,$$
for $x$ and $\beta$. Eliminating $x$, we obtain
$$3\tau(h^2)^2=2\tau(h^4)~,$$
which in turn yields
$$(1+\cos 2\pi\al)\beta^2=0~.$$
Since $\al$ is irrational this is valid only for $\beta=0$. Note that
$\beta=0$ does indeed correspond to an element with a connected spectrum,
namely $u+u^*$, which is the image of $v+v^*$ under a ``Fourier
transform''.\qed

We will show now that the eigenvectors of the operators $h(\del)$ depend
continuously on the spectral parameter $z$ in a sense to be made precise
below.

\bigskip
\noindent{\bf 2.14 Proposition.} Let $d>1$, and for every $m\in\N$ let
$|\del_m|\ge
d$, $z_m\in Sp(h(\del_m))$, $\xi^{(m)}\in\ell^2(\Z)\bsl\{0\}$ such that
$$h(\del_m)\xi^{(m)}=z_m\xi^{(m)}~,\quad\limsup_{|n|\to\infty}
\left|\xi^{(m)}_n\right|^{1/|n|}\le|\del_m|^{-1}\beta^{-1}$$
and $\lim\limits_{m\to\infty}z_m=z$.

Then there exist $\del\in\C$, $\xi\in\ell^2(\Z)\bsl\{0\}$ and $j_m\in\N$,
$c_m\in\C\bsl\{0\}$ such that $h(\del)\xi=z\xi$,
$\lim\limits_{m\to\infty}\lam^{2{j_m}}\del_m=\del$ and
$$\lim_{m\to\infty}\left\|D_\gam\left(c_m u^{j_m}
\xi^{(m)}-\xi\right)\right\|_\infty=0~~
\hbox{uniformly for}~d^{-1}\le\gam\le d~.$$

\bigskip
\Proof: It follows from (2.2) that there exists $t_m\in\C$ such that
$$kD_{\del_m}\xi^{(m)}=t_m D^{-1}_{\del_m}\xi^{(m)}~.$$
Let
$$\eta^{(m)}=\|D_{\del_m}\xi^{(m)}\|^{-1}_\infty D_{\del_m}\xi^{(m)}~.$$
Then $\|\eta^{(m)}\|_\infty=1$ and there exist $j_m\in\N$ such that
$$\left|\eta^{(m)}_{j_m}\right|\ge{1\over 2}~.$$
Adjusting $\eta^{(m)}$ suitably through shifting and multiplication by
scalars, we may assume that
$$\eta^{(m)}_0\ge{1\over 2}~.$$
Finally, switching to a subsequence if necessary, we may assume that the
sequence\hfill\break $\eta^{(1)},\eta^{(2)},\ldots$ converges weakly to
some $\eta\in\ell^\infty(\Z)\bsl\{0\}$. Adjusting the $t_m$ accordingly
we have
$$k\eta^{(m)}=t_m D^{-2}_{\del_m}\eta^{(m)}~.$$
Since $k\eta\not=0$, it follows that $\inf\{|t_m|\mid m\in\N\}>0$. Also,
since Proposition 2.12 implies that $|\del_1|,|\del_2|,\ldots$ is
convergent, $\sup\{|t_m|\mid m\in\N\}<\infty$. For, if this were not the
case, then we could find some point of density $\del$ for the $\del_m$
such that $D^{-2}_\del\eta=0$.

Switching once again to a subsequence if necessary, we may assume that
the sequences $\{\del_m\}$ and $\{t_m\}$ are convergent. Let
$\del=\lim\limits_{m\to\infty}\del_m$. Since $k$ defines a bounded
operator on $\ell^\infty(\Z)$, it follows that the corresponding sequence
of the $\xi^{(m)}$, after having been suitably scaled and adjusted
through shifts, does indeed converge uniformly to some eigenvector $\xi$
of $h(\del)$ for $z$, as claimed in the proposition. Since eigenvectors
of $h(\del)$ are essentially unique by Proposition 2.10, the conclusion
of the proposition applies to the original sequence of the $\xi^{(m)}$ as
well.\qed

\bigskip
\noindent{\bf Remark:} Tracing the steps in the proof of Lemmas 2.1 and
2.2 one can improve the convergence of eigenvectors in Proposition 2.14
to the effect that $\gam$ may range over any closed interval contained in
$(d^{-1}\beta^{-1},d\beta)$.\qed

Our next objective is to refine Propositions 2.9 and 2.13 by showing that
there is a natural parametrization of the eigenvalue problem for the
operators $h(\del)$, $|\del|>1$, through a Riemann surface $\tr$ covering
the resolvent set $\R=\C\bsl Sp(h)$. By doing so we will discover that
all the eigenvectors of these operators correspond to orbits of a cyclic
group of covering transformations on $\tr$, while the square of their
components can be obtained through evaluation of a single analytic
function on $\tr$. For the basic concepts of Riemann surfaces we are
going to use and for the terminology that comes with it, we refer to
[AS], Chapters I and II.

Let $I$ be the smallest interval containing $\sph$. There exists a unique
analytic function $G$ on $\C\bsl I$ such that
$$\log|G(z)|=\int\log|z-s|d\mu(s)-\log\beta~,
\quad G(\r^{+}\bsl\sph)\subset\r^{+}~.$$
If $f$ is a closed arc surrounding a component $K$ of $\sph$, then any
two analytic continuations of $G$ along $f$ over the same point $z$ on
$f$ differ by a multiplicative constant of the form
$$e^{2\pi\mu(K)ni}~.$$
A well-known fact in the $K$-theory of the irrational rotation
$C^*$-algebra provides us with the information that
$$\mu(K)\in\(\Z+\al\Z\)\cap[0,1]~,$$
so that the said multiplicative factor takes the form $\lam^{2n}$. Now
$G$ can be continued analytically along any arc in ${\R}={\C}\bsl{\sph}$.
So, if we let $\F$ be the universal covering of $\R$, $\G$ the
corresponding group of covering transformations, and finally $\tg$ the
analytic function on $\F$ obtained through analytic continuation of $G$,
then we have for every $g\in\G$ 
$$\tg\circ g=\lam^{2n}\tg~~\hbox{for some}~n\in\Z~.$$
Let
$$\G_0=\{g\in\G\mid\tg\circ g=\tg\}~.$$
Then $\G_0$ is a normal subgroup of $\G$ with a cyclic quotient group
$\G/\G_0$. The $\G_0$-orbits in $\F$ form a Riemann surface $\tr$ which
covers $\R$, and whose group of covering transformations corresponds to
$\G/\G_0$ in a natural way. Let $p$ be the covering map of $\tr$ over
$\R$. Finally, $\tg$ defines an analytic map on $\tr$ which we denote by
$G$ again. We assume $\om$ to be that generator of the group of covering
transformations of $\tr$ over $\R$ for which $G(\om(z))=\lam^2 G(z)$.

We are now in a position to state the major claim in this paragraph. For
its proof, we need another technical lemma.

\bigskip
\noindent{\bf 2.15 Lemma.} Suppose $f:\tr\to\T$ has the following
properties
\item{(i)} Every subsequence of $\{f^n\mid n\in\Z\}$ has a
subsequence which converges uniformly on compact subsets of $\tr$.
\item{(ii)} If $\lim\limits_{n\to\infty}z_n=z$ and
$\lim\limits_{n\to\infty}f(z_n)=c$, then $c\in\{f(z),\lam^2f(z)\}$.

\noindent Then $f(\tr)\subset\{x\lam^{2n}\mid n\in\Z\}$ for some
$x\in\T$.

\bigskip
\Proof: Let $K\subset\tr$ be compact and connected. It suffices to show
that $f(K)\subset\{x\lam^{2n}\mid n\in\Z\}$ for some $x\in\T$. Let $\Om$
be the uniform closure in $\ell^\infty(K)$ of $\{g^n\mid n\in\Z\}$, where
$g=f/k$. It follows from (i) that $\Om$ is a metrizable and compact
group. To every point $y$ in $g(K)$ there corresponds a character $\vp_y$
of $\Om$, a continuous homomorphism from $\Om$ into $\T$, such that
$\vp_y(g)=y$. Let $\Y=\{\vp_y\mid y\in g(K)\}$. Since the uniform
structure on $\Om$ is inherited from $\ell^\infty(K)$, the set $\Y$ is
equicontinuous on $\Om$. Hence $\oY$, the closure of $\Y$ in the topology
of uniform convergence on $\Om$, is compact. In other words, $\oY$ is a
compact subset of the dual group $\widehat\Om$ of $\Om$. Since $\Om$ is
compact,$\widehat\Om$ is discrete. Therefore $\oY$, being a compact
subset of $\widehat\Om$, must be finite. It follows that $g(K)$ is a
finite subset of $\T$.

For $x\in\T$, let $\calo_x=\{x\lam^{2n}\mid n\in\Z\}$. Now let $x\in\T$
be such that
$$M=g(K)\cap\calo_x\not=\emptyset~.$$
Let $N=g(K)\bsl\ \calo_x$. Since $g(K)$ is finite, $M$ and $N$ are finite
as well and hence, closed. Therefore, since $K$ is compact, property (ii)
implies
$$\ovl{g^{-1}(M)}\subset M\cup\lam^2 M~,\quad \ovl{g^{-1}(N)}
\subset N\cup\lam^2N~.$$
By the definition of $M$ and $N$
$$(M\cup\lam^2M)\cap(N\cup\lam^2 N)=\emptyset~.$$
Hence,
$$\ovl{g^{-1}(M)}\cap\ovl{g^{-1}(N)}=\emptyset~,\quad \ovl{g^{-1}(M)}
\cup\ovl{g^{-1}(N)}=K~.$$
Since $K$ is connected by assumption,
$$g^{-1}(M)=\ovl{g^{-1}(M)}=K~,$$
which settles the claim.\qed

\bigskip
\noindent{\bf 2.16  Theorem.} For every $z\in\tr$, there exists
$\xi^{(z)}\in\ell^2(\Z)\bsl\{0\}$
such that
$$h(G(z))\xi^{(z)}=p(z)\xi^{(z)}~,\quad {\limsup_{|n|\to\infty}}
\left|\xi_n^{(z)}\right|^{1/|n|}
\le|G(z)|^{-1}\beta^{-1}~,$$
and for every $p\in\Z$, $\vartheta_p(z)=\xi^{(z)}_0\xi^{(z)}_p$ is an
analytic function on $\tr$ with the property
$$\vartheta_p(\om^n(z))=\xi^{(z)}_n\xi^{(z)}_{n+p}~.$$

\bigskip
\Proof: Propositions 2.9 and 2.12 taken together imply that for every
$z\in\tr$, there exist $x_z\in\T$, $\xi^{(z)}\in\ell^2(\Z)\bsl\{0\}$ such
that
$$h(x_zG(z))\xi^{(z)}=p(z)\xi^{(z)}~,\quad \ovl{\lim_{|n|\to\infty}}
\left|\xi^{(z)}_n\right|^{1/|n|}\le|G(z)|^{-1}\beta^{-1}~.$$
On account of (2.2), there exists for every $z\in\tr$ a $t_z\in\C$ such
that
$$kD_{x_zG(z)}\xi^{(z)}=t_z D^{-1}_{x_zG(z)}\xi^{(z)}~.\leqno(*)$$

Applying $u^m$ to both sides of this identity yields
$$kD_{x_zG(z)}u^m\xi^{(z)}=(x_zG(z))^2t_z D^{-1}_{x_zG(z)}u^m\xi^{(z)}~.$$
Therefore, we can adjust $x_z$, $\xi^{(z)}$ and $t_z$ in $(*)$ such that
$$1\le|t_z|<|G(z)|^2~.$$
It follows from Proposition 2.10 that
$$\lim_{n\to\infty}z_n=z~,~~\lim_{n\to\infty}|t_{z_n}|
=s~~\hbox{and}~~\lim_{n\to\infty}x_{z_n}=y$$
implies
$$s\in\left\{|t_z|, |G(z)|^2|t_z|\right\}~,\quad
y\in\left\{x_z,\lam^2 x_z\right\}~.$$ Now let
$$\vp^{(z)}(a)=\langle  a\xi^{(z)},\ovl{\xi^{(z)}}\rangle~,\quad a\in\A~.$$
Then
$$\vp^{(z)}(h(x_zG(z))a)=\vp(ah(x_zG(z))=p(z)\vp^{(z)}~.$$
So, if we define
$$\psi^{(z)}_{pq}=\ovl{x}^q_zG(z)^{-q}\vp(w_{pq})~;\quad p,q\in\Z~,$$
then the double sequence $\left\{\psi^{(z)}_{pq}\right\}$ is a solution
to the system of difference equations in the Appendix A2. Also,
$\left\{\psi^{(z)}_{pq}\right\}$ decays exponentially uniformly in $p$
with an order less than or equal to $\beta^{-1}$. Therefore, it follows
from A2 that, after scaling the eigenvectors $\xi^{(z)}$ suitably,
$$\lam^{pq}\ovl{x}^q_z\sum^\infty_{n=-\infty}\ovl{\lam}^{2qn}
\xi^{(z)}_n\xi^{(z)}_{n+p}
=\left(c_{pq}(p(z))-d_{pq}(p(z)\right)G(z)^q~;\quad p,q\in\Z~,\leqno(**)$$
where $c_{pq}$ and $d_{pq}$ are analytic functions on $\R$ which are
determined by the identities
$$\eqalign{h-x)^{-1}&=\sum^{\infty}_{p,q=-\infty}c_{pq}(x)w_{pq}
\quad\hbox{and}\cr
\noalign{\vskip8pt}
(h(\del)-x)^{-1}&=\sum^{\infty}_{p,q=-\infty}d_{pq}(x)\del^qw_{pq}~,\quad
|\del|>|G(z)|~.}$$

If $d>1$ and $K\subset\{z\in\tr\mid |G(z)|\ge d\}$ is compact, then
Proposition 2.11 ensures that
$$\|D_\gam\xi^{(z)}\|_\infty~~\hbox{is uniformly bounded in~}
~z\in K, d^{-1}\le|\gam|\le d~.$$ Since the functions on the right hand
side of (2.4) are analytic on $\R$ as well as uniformly bounded on $K$ in
$p$ and $q$, we conclude that if we let $f(z)=\ovl{x}_z$, then any
subsequence of $\{f^q\}$ has a subsequence which is locally uniformly
convergent. Thus we have show that $f$ has the properties (i) and (ii) in
Lemma 2.15. It follows that
$$f(\tr)\subset\{\oy\lam^{2n}\mid n\in\Z\}~~\hbox{for some}~y\in\T~.$$
By applying suitable powers of $u$ to the vectors $\xi^{(z)}$, we can
adjust $(**)$ to the effect that $x_z=y$ for every $z\in\tr$.

Let $t>\max Sp(h)$. Then there exists $z_t\in\tr$ such that
$G(z)\in\r^{+}$. Hence
$$h(yG(z_t))\xi^{(z_t)}=t\xi^{(z_t)}~,\quad h(\oy G(z_t)\J\ovl{\xi^{(z_t)}}
=t\J\ovl{\xi^{(z_t)}}$$
and by Proposition 2.10, $\oy=\lam^{2m}y$ for some $m\in\Z$. Thus,
adjusting the eigenvectors through a suitable power of the shift we may
assume that $y\in\{1,-1,\lam,-\lam\}$. In order to show that $y$ equals
actually $1$, one can proceed as follows: Since
$\beta^{-1}G(z_t)h(yG(z_t))$ approaches $yv$ as $t\to\infty$, one can use
$(**)$ to show that $\xi^{(z_t)}$ approaches an eigenvector $\eta$ of
$yv$ (in $\ell^2(\Z)$, say) for some positive eigenvalue $s$. Thus
$s\in\{\oy\lam^{2n}\mid n\in\Z\}$. The only way this can happen is when
$y=1$.

Proposition 2.13 in conjunction with $(**)$ entails that the function
$$\psi_p:\T\times\tr\to\C~;\quad \psi_p(x,z)
=\sum^{\infty}_{n=-\infty}x^n\xi^{(z)}_n\xi^{(z)}_{n+p}$$
is continuous and analytic in $z$ for every $x\in\T$. It follows that the
functions
$$\vt_p(z)=\int\psi_p(x,z)\,dx=\xi^{(z)}_0\xi^{(z)}_p$$
are analytic. Moreover, since $G(\om(z))=\lam^2G(z)$,
$$\vt(\om^n(z))=\xi_n^{(z)}\xi^{(z)}_{n+p}$$
by $(**)$, as claimed.\qed

\bigskip
\noindent{\bf Remarks.} 1) The second part of Theorem 2.16 says (take
$p=0$), that given $\xi^{(z)}_0$ in a neighborhood of $z$, one can
generate $\(\xi^{(z)}_n\)^2$ through analytic continuation of
$\(\xi^{(z)}_0\)^2$.
\bigskip
2) Theorem 2.16 also shows that we can identify the points of the Riemann
surface $\tr$ with the one-dimensional eigenspaces of the operators
$h(G(z))$ for $p(z)$ in such a way that the covering transformation $\om$
corresponds to the two-sided shift on $\ell^2(\Z)$.
\bigskip
3) From $(**)$ in the proof of Theorem 2.16 we can extract the identity
$$\sum^{\infty}_{n=-\infty}\(\xi^{(z)}_n\)^2=\int(t-p(z))^{-1}d\mu(t)~.$$
It follows that $\sum\limits^{\infty}_{n=-\infty}\(\xi^{(z)}_n\)^2$
equals zero if and only if $z$ is a critical point for the conductor
potential of $Sp(h)$. These critical points are simple; and exactly one
is located in every gap of $Sp(h)$ and nowhere else. This shows that the
eigenvectors $\xi^{(z)}$ carry the relevant information regarding the gap
structure for $Sp(h)$ in a very explicit form.
\bigskip
4) Elaborating on the comments made above, there are some implications
for the case $|\del|=1$. One might be tempted to try to generate
exponentially decaying eigenvectors for the operators $h(x)$, $|x|=1$, as
(uniform) limits of the eigenvectors $\xi^{(z)}$. Assuming that $Sp(h)$
has infinitely many gaps, one is confronted with the following
impediment: Let $z_1,z_2,\ldots$ be a sequence of critical points
converging to $\chi\in Sp(h)$. Suppose that $h(x)\xi=\chi\xi$ for some
$\xi\in\ell^2(\Z)\bsl\{0\}$. Then
$$\sum^{\infty}_{n=-\infty}\xi^2_n=y\|\xi\|_2\not=0
\quad\hbox{for some}~y\in\T~.$$
It follows that $\xi$ cannot be approximated by the (suitably scaled)
sequence $\xi^{(z_1)},\xi^{(z_2)},\ldots$ in the Hilbert space norm.

\bigskip
\noindent{\bf 2.17 Corollary.} There exist homeomorphisms $\sig$ and $\iota$
on $\tr$ such that
$$\eqalign{\xi^{(\sig(z))}&=\J\ovl{\xi^{(z)}}\quad~,~p(\sig(z))=\ovl{p(z)}~,
~G(\sig(z))=\ovl{G(z)}\cr
\xi^{(\iota(z))}&=D_{-1}\xi^{(z)}~,~p(\iota(z))=-p(z)~,~G(\iota(z))=-G(z)}~.$$
In particular, $\sig^2=\iota^2=Id$, $\sig\circ\om=\om^{-1}\circ\sig$,
$\om\circ\iota=\iota\circ\om$.

\bigskip
\Proof:  This follows immediately from Theorem 2.16, noting that
$$h(\ovl{G(z)})\J\ovl{\xi^{(z)}}=\ovl{p(z)}\J\ovl{\xi^{(z)}}$$
and
$$h\(-{G(z)}\)D_{-1}\xi=-\ovl{p}(z)D_{-1}\xi^{(z)}~.\qeddis

The following is in preparation for the discussion in paragraph 3.

\bigskip

\noindent{\bf 2.18 Corollary.} There exists an analytic function $\Gam$
on $\wt\R$ which has the following properties:
$$kD_{G(z)}\xi^{(z)}=\Gam(z)D_{G(z)^{-1}}\xi^{(z)},~\Gam(\om(z))=G(z)^2\Gam(z),
~\Gam(\sig(z))=\ovl{\Gam(z)}^{-1}$$
$$~\hbox{and}~\Gam(\iota(z))=\Gam(z).$$

\Proof: There clearly exists a function satisfying the stated
identities. To show that it is actually analytic, we consider scalar
products
$$\langle kD_{G(z)}\xi^{(z)}, a\ovl{\xi^{(z)}}\rangle=\Gam(z)
\langle D_{G(z)^{-1}}\xi^{(z)}, a\ovl{\xi^{(z)}}\rangle~,$$
where $a$ is any linear combination of the elements $w_{pq}$. Written out
in components, one can see that the function on the left as well as the
second factor on the right are algebraically generated by $G$, $G^{-1}$
and functions of the form $\vt_p\circ\om^n$, all of which are analytic.
Therefore, $\Gam$ has to be analytic as well.\qed

\bigskip

\item{\bf 3.}{\bf The kernel of a family of related operators}
\medskip
The first identity in Corollary 2.18 gives rise to an eigenvalue problem
in its own right, involving the operator $k$ and the function $\Gam$ as
an eigenvalue parameter. The question arises whether the two related
eigenvalue problems are actually equivalent. It will be shown that this
is essentially the case. More specifically, it will be shown that to any
two distinct points $z_1$ and $z_2$ in the complement of a (possibly
empty) discrete subset $\wt\R$, there correspond distinct pairs of
parameters $(G(z_1),\Gam(z_1))$ and $(G(z_2), \Gam(z_2))$. To this end,
we introduce new operators whose kernel will hold all the relevant
information.

Let $g$ be the function defined in (1.6) and let $\gam\in\r$ be close to
$1$ having the property that
$$g(\gam\T)\cap\Bigl(\Bigl\{\pi\Bigl(n+{1\over 2}\Bigr)
~\Bigl|~~ n\in\Z\Bigr\}\Bigr)
\cup\left\{\pi\(m+i\(n+{\pi\over 4}\)\)~\Bigl|~~ m,n\in\Z\right\}=\emptyset~.$$
Let $\Gam$ be as in Corollary 2.18 and let
$$\Om=\{z\in\tr \mid 1+\Gam(\om^n(z))\not=0~~\hbox{for all}~n\in\Z\}~.$$
Since $\Gam$ is analytic, $\Om$ is a discrete subset of $\tr$. In the
following, $\Om$ will be augmented as needed, but it will always be
discrete. Now we define for every $z\in\tr\bsl\Om$ a bounded operator
$H_\gam(z)$ on the Hilbert space $\ell^2(\Z)$,
$$\(H_\gam(z)\xi\)_n=\sum^{\infty}_{j=-\infty}a_j\xi_{n+j}
+i{1-\Gam(\om^n(z))\over{1+\Gam(\om^n(z))}}\xi_n~,$$
where
$$\tan g(\gam u)=\sum^{\infty}_{j=-\infty}a_j u^j~.$$
The sequence $\{a_j\}$ decays exponentially as $|j|\to\infty$. Obviously,
$H(z)$ is analytic in $z$.

Moreover,
$$\left\{\eqalign{H_\gam(\om(z))&= uH_\gam(z)u^*\cr
H_\gam(\sig(z))&=\J H_{\gam^{-1}}(z)^*\J~.}\right.\leqno(3.1)$$

\noindent{\bf 3.1 Proposition.} The operator $H(z)$ is Fredholm with
index zero. Any two of these operators are compact perturbations of each
other, and their essential spectrum equals
$$(\tan g(\gam\T)+i)\cup(\tan g(\gam\T)-i)~.$$

\Proof: Let $T\in\B(\ell^2(\Z))$ be defined as follows:
$$(T\xi)_n=\sum^{\infty}_{j=-\infty}a_j\xi_{n+j}+\E(n)\xi_n~,$$
where
$$\E(n)=\cases{-i~~,&$n\ge 0$\cr
\noalign{\vskip4pt}
~~i~~,&$n<0$.}$$
\medskip
\noindent Since $|\Gam(\om^n(z))|=|G(z)^{2n}\Gam(z)|$ and $|G(z)|>1$,
\smallskip
$$\lim_{n\to\infty}i{1-\Gam(\om^n(z))\over{1+\Gam(\om^n(z))}}=-i~~,~~
\lim_{n\to-\infty}i{1-\Gam(\om^n(z))\over{1+\Gam(\om^n(z))}}=i~.$$
\smallskip
\noindent It follows that $H(z)$ is a compact perturbation of $T$. Let
$$T_0=T-\E~,$$
and let $\tt_0$, $\te$ be the range of $T_0$, $\E$, respectively, in the
Calkin algebra $\B(\ell^2(\Z))/\K$, where $\K$ denotes the $C^*$-algebra
of compact operators on $\ell^2(\Z)$. Then $\tt_0$ and $\te$ are normal
operators which commute and their joint spectrum equals
$$\tan g(\T)\times\{-i,i\}~.$$
It follows that
$$Sp\(\tt_0+\te\)=(\tan g(\gam\T)+i)\cup(\tan g(\gam\T)-i)$$
which equals the essential spectrum of $T$. By our choice of $\gam$,
$\tt_0+\te$ is invertible, which entails that $H(z)$ is Fredholm with
index zero.\qed

\bigskip
In the following we will constantly make use of Theorem 2.16. In
particular, the eigenvectors which occur will be assumed to decay of a
sufficiently high order, so that all manipulations make sense, unless
specified otherwise. The significance of the operators $H(z)$ rests with
the following statement.

\bigskip
\noindent{\bf 3.2 Proposition.} Suppose $h_\gam(G(z))\xi=p(z)\xi$ for
some $z\in\tr\bsl\Om$. Then there exist $t_n\in\C$ such that
$$t^2_n=\Gam(\om^n(z))~,~~H_\gam(z)\eta=0~,$$
where
$$\eta_n=\(t_n+t^{-1}_n\)\xi_n~.$$

\Proof: By Corollary 2.18
$$k_\gam D_{G(z)}\xi=\Gam(z)D_{G(z)^{-1}}\xi~.$$

Rearranging the diagonal operators involved, we obtain
$$Te^{ig(\gam u)}\xi=T^{-1}\xi$$
where $T$ is a diagonal operator such that $(T\xi)_n=t_n\xi_n$ with
$t^2_n=\Gam(\om^n(z))$. Using the identity
$$e^{ig(\gam u)}={1+i\tan{g(\gam u)\over{2}}\over{1-i\tan{g(\gam u)\over{2}}}}$$
we have
$$\(1+i\tan{g(\gam u)\over{2}}\)T\xi=\(1-i\tan{g(\gam u)\over{2}}\)T^{-1}\xi$$
or
$$t_n\xi_n+i\sum^{\infty}_{j=-\infty}a_jt_{n+j}\xi_{n+j}
=t^{-1}_n\xi_n-i\sum^{\infty}_{j=-\infty}a_jt^{-1}_{n+j}\xi_{n+j}$$
which turns into
$$\sum^{\infty}_{j=-\infty}a_j\eta_{n+j}
+i{t_n^{-1}-t_n\over{t_n^{-1}+t_n}}\eta_n=0~.$$ Since $t^2_n=\Gam(z)$,
the second term equals
$$i{1-\Gam(\om^n(z))\over{1+\Gam(\om^n(z))}}~.\qeddis

\noindent{\bf Remark:} If we replace $\xi$ by $\txi=D_{-1}\xi$ in
Proposition 3.2, then
$$h(G(\iota(z))\txi=-p(z)\txi~.$$
But if we transform $\txi$ analogous to $\xi$, thus obtaining a vector
$\teta$, we have $\eta=\teta$. Succinctly put, one can say that to every
non-zero element in the kernel of $H_\gam(z)$ there correspond
eigenvectors of $h(G(z))$ and $h(G(\iota(z))$, respectively, whose
eigenvalues differ by a negative sign. Also note that $H(\iota(z))=H(z)$,
since $\Gam(\iota(z))=\Gam(z)$. This means that the operators $H(z)$ are
more appropriately parametrized through the Riemann surface obtained from
$\tr$ by identifying $z$ and $\iota(z)$.\qed

Propositions 3.1 and 3.2 taken together say that $H(z)$ is a Fredholm
operator of index zero with a non-trivial kernel and that $0$ is an
isolated point in the spectrum of $H_\gam(z)$ for every $z\in\tr\bsl\Om$.
This situation provides the proper setting for the employment of analytic
perturbation theory of linear operators as expounded in [Ko], for
instance.

First we enlarge $\Om$ by the set of points $z\in\tr$ for which the
following is true: For every neighborhood $\U\subset\tr$ of $z$ and for
every neighborhood $\V\subset\C$ of $0$, there exists $\tz\in\U$ such
that $H_\gam(\tz)$ has a non-zero eigenvalue in $\V$. Since $H_\gam(z)$
is analytic in $z$, the inclusion of those branch-points in $\Om$ still
yields a discrete subset of $\tr$.

For every point $z\in\tr\bsl\Om$, let
$$P(z)={1\over{2\pi i}}\oint\(t-H_\gam(z)\)^{-1}dt~,$$
where the integral is taken over a positively-oriented circle enclosing
$0$ but no other point in the spectrum of $H_\gam(z)$. Then $P(z)$ is a
projection of finite rank. Moreover, $P(z)$ is analytic in $z$ and the
rank of $P(z)$ is constant. For every $z\in\tr\bsl\Om$, let
$$N(z)=H_\gam(z)P(z)~.$$
Then $N(z)$ is a nilpotent operator which is also analytic in $z$.

Our next goal is to show that the kernel of $H_\gam(z)$ is
one-dimensional for every\break $z\in\tr\bsl\Om$. In preparation of this,
we settle a number of technical questions first.
\bigskip
\noindent{\bf 3.3 Lemma.} There exists a discrete subset $\Om_0\subset\tr$
such that the following holds true: If $z\in\tr\bsl\Om_0$ and
$\L\subset\ell^2(\Z)$ is a finite-dimensional subspace which is invariant
under $h(G(z))$, then $\L$ contains a linear basis of eigenvectors of
$h(G(z))$.

\bigskip
\Proof: Let $\xi^{(z)}$ and $\vt_p$ be as in Theorem 2.16. Let
$$\Om_0=\{z\in\tr\mid\hbox{There exists $\tz\in G^{-1}(\{p(z)\})$
such that $\vt_0(\tz)=0$}\}~.$$ Since $\vt_0$ is analytic, the set
$M_0=\{\tz\mid\vt_0(\tz)=0\}$ is a discrete subset of $\tr$. For every
compact subset $K\subset\tr$, the set $\{z\in K\mid z\in
G^{-1}(\{p(\tz)\})$for some $\tz\in M_0\}$ is finite. Hence $\Om_0$ is
discrete.

Now let $z\in\tr\bsl\Om_0$ and consider the Jordan canonical form of
$h(G(z))$ restricted on $\L$. We need to show that the nilpotent
components in this decomposition are trivial. To this end,we need to show
that if $\tz\in\tr$ such that $p(\tz)$ is an eigenvalue for $h(G(z))$
with eigenvector $\xi^{(\tz)}$, then there does not exist
$\eta\in\ell^2(\Z)$ such that $(h(G(z))-p(\tz))\eta=\xi^{(\tz)}$. Suppose
that the opposite is true: There exists $\eta$ with the said property.
Then
$$\eqalign{\vt_0(z)
&=\langle\xi^{(\tz)},\ovl{\xi^{(\tz)}}\rangle
=\langle\(h(G(z))-p(\tz)\)\eta,\ovl{\xi^{(\tz)}}\rangle\cr
&=\langle\eta,\(h(G(z))-p(\tz)\)^*\ovl{\xi^{(\tz)}}\rangle\cr
&=\langle\eta,\ovl{\(h(G(z))-p(\tz)\)\xi^{(\tz)}}\rangle=0~,}$$ which
means that $z\in\Om_0$, thus contradicting our assumption on $z$.\qed

We augment $\Om$, if necessary, to include the set $\Om_0$.

\medskip
\noindent{\bf 3.4 Corollary.} Let $z\in\tr\bsl\Om$ and for every
$n\in\Z$, let $t_n$ be chosen as in Proposition 3.2. Let
$$\L=\{\xi\mid\xi_n=(t_n+t^{-1}_n)^{-1}\eta_n~,
~\hbox{where}~\eta\in\ell^2(\Z), H_\gam(z)\eta=0\}~.$$ Then $\L$ contains
a linear basis consisting of eigenvectors of $h_\gam(G(z))$.

\bigskip
\Proof: Since the kernel of $H_\gam(z)$ is finite-dimensional,
in view of Lemma 3.3 all that needs to be shown is that
$h_\gam(G(z))\L\subset\L$. Let $\eta\in\ell^2(\Z)$, $H_\gam(z)\eta=0$ and
let $\xi_n=(t_n+t^{-1}_n)\eta$, $\txi=h_\gam(G(z))\xi$,
$\teta_n=(t_n+t^{-1}_n)\txi_n$. Since $|t_n|\,|G_n(z)|^{-n}$ is constant,

$$\sum^{\infty}_{n=-\infty}|t_n+t^{-1}_n|^2|\xi_n|^2<\infty~\hbox{implies}~
\sum^{\infty}_{n=-\infty}|t_n+t^{-1}_n|^2|\txi_n|^2<\infty~,$$
which implies $\teta\in\ell^2(\Z)$.

Next, since $D_{G(z)}k_\gam D_{G(z)}$ and $h_\gam(G(z))$ commute, and
since $D_{G(z)}k_\gam D_{G(z)}\xi=\Gam(z)\xi$,
$$D_{G(z)}k_\gam D_{G(z)}\txi=\Gam(z)\txi~,$$
which in turn implies that $H_\gam(z)\teta=0$. In conclusion,
$\txi=h_\gam(G(z))\xi\in\L$.\qed

Since $N$ is analytic, the rank of $N(z)$ is constant on $\tr\bsl\Om$
with the possible exception of a discrete subset. We now enlarge $\Om$ by
this set of exceptional points.

\bigskip
\noindent{\bf 3.5 Lemma.} For every $z\in\tr\bsl\Om$, let $Q(z)$ be the
projection whose kernel equals the range of $H_\gam(z)$ and whose range
equals the kernel of $H_\gam(z)$. Then $Q$ is analytic.

\bigskip
\Proof: Let $z_0\in\tr\bsl\Om$ and let $\U\subset\tr\bsl\Om$ be a simply
connected open neighborhood of $z_0$. We shall shrink $\U$ when needed.
By [Ko], II--\S4.2 and VII--\S1.3, there exists an analytic function $T$
from $\U$ into $\B(\H)$ such that $T(z)$ is invertible for every $z\in\U$
and
$$T(z)P(z_0)T(z)^{-1}=P(z)~.$$
Let
$$\tn(z)=T^{-1}(z)N(z)T(z)~.$$
Then $\tn$ is analytic on $\U$ and $\tn(z)$ is nilpotent.

We may consider $\tn(z)$ as an operator on the finite-dimensional
subspace $\L=P(z_0)\ell^2(\Z)$. We choose a linear basis in $\L$, and we
denote the matrix representation of $\tn(z)$ with respect to this basis
once again by $\tn(z)$. Let $m$ be the dimension of $\L$ and let $\B$ be
the canonical basis of $\C^m$. Then we choose a subset $\B_1$ of $\B$
such that $\{\tn(z_0)e\mid e\in\B\}$ is a linear basis for
$\tn(z_0)\C^m$. Since the functions $z\mapsto\tn(z)e$ are analytic in
$\U$ for all $e\in\B_1$, $\B_1(z)=\{\tn(z)e\mid e\in\B_1\}$ is a basis
for $\tn(z)\C^m$ for every $z$ in an open neighborhood of $z_0$ contained
in $\U$. We replace $\U$ by that neighborhood. Since $\tn$ is nothing but
a matrix whose entries are analytic functions, we can choose a submatrix
$M$ of $\tn$ of type $r\times m$, where $r$ equals the rank of
$\tn(z_0)$, such that $M(z_0)$ has rank $r$. Again, since the entries of
$M$ are analytic in $z$, $r=\rank(M(z))=\rank(\tn(z))$ for every $z$ in
an open neighborhood of $z_0$ contained in $\U$. Once again we replace
$\U$ by that neighborhood. We can now use Cramer's rule to solve a system
of linear equations for analytic functions $f_1,\ldots,f_r$ on $\U$ with
values in $\C^m$ such that $\B_2(z)=\{f_1(z),\ldots,f_r(z)\}$ is a basis
for the kernel of $\tn(z)$. Let $\tq(z)$ be the projection whose range
equals the kernel of $\tn(z)$ and whose kernel equals the range of
$\tn(z)$. It follows that $\tq$ is analytic on $\U$. By construction
$Q(z)=T(z)(\tq(z)\circ P(z_0))T(z)^{-1}$. In conclusion, $Q$ is locally
analytic and hence analytic on $\tr\bsl\Om$ as claimed.\qed

\bigskip
\noindent{\bf 3.6 Lemma.} Let $F$ be a mapping from $\R\bsl p(\Om)$ into
the subsets of $\R$ such that $F(z)$ contains exactly $r$ elements.
Suppose that $F$ has the following properties:
\item{(I)} If $z\in Sp(h(\del))$, then $F(z)\subset Sp(h(\del))$.
\item{(II)} For every $z_0\in\R$, there exists an open neighborhood $\U$
of $z_0$ and for every $j$, $1\le j\le r$, there exists an analytic
function $s^{(z_0)}_j$ on $\U$ such that
$$F(z)=\{s_1^{(z_0)}(z),\ldots,s^{(z_0)}_r(z)\}~~\hbox{for every $z\in\U$}~.$$

\noindent Then $F(z)\subset\{z,{-}z\}$. (A mapping $F$ having the property (II) is
called an {\it analytic multifunction\/}).

\bigskip
\Proof: For every $z\in\R\bsl p(\Om)$, we form the polynomial
$$\(X-s_1^{(z)}(z)\)\cdot\ldots\cdot\(X-s_r^{(z)}(z)\)
=X^r+\sum^{r-1}_{j=0}f_j(z)X^j~.$$
Then (II) entails that the coefficients $f_j$ are analytic functions on
$\R\bsl p(\Om)$. By (I) the sets $F(z)$ are uniformly bounded in a
punctured neighborhood of any point in $p(\Om)$. It follows that the same
is true for the functions $f_j$. Thus, all points in $p(\Om)$ are
removable singularities for $f_j$, and we may consider $f_j$ as an
analytic function on $\R$. It follows that the roots of the above
polynomial are branches of analytic functions with exceptional point
located in $p(\Om)$ (see [Ko], II--\S1.2).

Let $d$ be the smallest number such that the logarithmic potential
associated with the equilibrium distribution $\mu$ of $Sp(h)$ has a
critical point $y$ on the corresponding level curve. More precisely,
$$\int\log|y-t|d\mu(t)=d~,$$
while
$$\int(y-t)^{-1}d\mu(t)=0.$$
By Corollary 2.13, we know that there always exists a critical point. Let
$$\R_d=\Bigl\{z\in\C~\Bigl|~\int\log|z-t|d\mu(t)>d\Bigr\}~.$$
Then there exists an analytic function $G_d$ on $\R_d$ such that
$$\log|G_d(z)|=\int\log|z-t|d\mu(t)$$
and
$$G_d(z)\subset\r^{+}~\hbox{whenever}~z\in\r^{+}~.$$
Also, $G_d$ has a simple pole at infinity. Let $\oR_d$ be the closure of
$\R_d$ and let $M$ be the set of critical points in $\oR_d$. Then $G_d$
has a continuous extension on $\oR_d\bsl M$ and it approaches two
distinct points as $z$ approaches a point in $M$. Let $z_0\in\R_d\bsl
p(\Om)$. Then (I) entails that
$$\left.\eqalign{{}&s_j^{(z_0)}(z)\subset\R_d~,\cr
{}&\left|G_d\(s_j^{(z_0)}(z)\)\right|=\bigl|G_d(z)\bigr|}\right\}~~\hbox{in
a neighborhood $U$ of $z_0$}$$ for $1\le j\le r$.

It follows that there exist $c_j\in\T$ such that
$$G_d\circ s^{(z_0)}_j=c_jG_d~\hbox{on}~\U~,\quad 1\le j\le r~.$$
Since this holds true for every $z_0\in\R_d\bsl p(\Om)$, (II) implies
that the set of those scaling factors $c_j$ is the same for all points in
$\R_d\bsl p(\Om)$. Hence
$$G_d\circ F(z)=\{c_1,\ldots,c_r\}G_d(z)~,$$
for $z\in\R_d\bsl p(\Om)$, where $G_d\circ F$ denotes the multifunction
obtained by applying $G_d$ to every element in $F(z)$. Since $G_d$ is a
conformal map from $\R_d\cup\{\infty\}$ onto \hfill\break $\{z~\bigl|~
|z|>e^d\}\cup\{\infty\}$, it follows that there exist analytic functions
$s_1,\ldots,s_r$ on $\R_d$ such that
$$F(z)=\{s_1(z),\ldots,s_r(z)\},\quad z\in\R_d\bsl p(\Om)~.$$
Every $s_j$ is transformed via $G_d$ into a rotation by $c_j$ on $\{z\mid
|z|>e^d\}$. Moreover, the analytic continuation of $s_j$ across the
boundary of $\oR_d$ does not result in the occurrence of exceptional
points (branching-point) in $\oR_d\bsl\R_d$ or points of order larger
than one. For, the existence of such points would conflict with the fact
that all the critical points in $M$ are simple. It follows that $s_j$ has
a differentiable extension on $\oR_d$ with non-vanishing derivatives. So,
$s_j$ has an inverse transformation with the same properties. Let $\K$ be
the group of transformations on $\oR_d$ which is generated by
$s_1,\ldots,s_r$ and their inverses. Since the derivative $G'_d(z)$
approaches zero as $z$ approaches a point in $M$, it follows that every
$s\in\K$ maps $M$ into $M$. Moreover, if $z\in M$ is a fixed point for
some transformation in $\K$, that transformation must be the identity.
Hence, since $M$ is finite, $\K$ is finite too. In particular, the group
of scaling factors
$$\cc=\{c\in\T~\big|~G_d\circ s=cG_d~\hbox{for some $s\in\K$}\}$$
is finite, which in turn means that $c^m=1$ for every $c\in\cc$, where
$m$ is the cardinality of $\cc$. We want to show that
$\cc\subset\{-1,1\}$.

To this end, we first note that $G_d(z)\subset G(p^{-1}(z))$, if we let
$G_d(z)$ be the set of the two limit points of $G_d$ at $z$ in case $z$
is in $M$. So, let $z\in M$. Since $z$ is real, since
$G\circ\sig=\ovl{G}$, and since $\sig\circ\om=\om^{-1}\circ\sig$, there
exists a point $\tz\in\P^{-1}(z)$ such that either
$$\sig(\tz)=\tz~~{\rm or}~~\sig(\tz)=\om(z)~.$$
In either case, $G(p^{-1}(z))\subset\{\lam^n\mid
n\in\Z\}\cup\{-\lam^n\mid n\in\Z\}$. Now let $s\in\K$ and let $z'=s(z)$.
Then
$$G_d(z), G_d(z')\subset\{\lam^n~\big|~n\in\Z)\cup\{-\lam^n~\big|~n\in\Z\}~.$$
Let $G_d\circ s=cG_d$ for some $c\in\cc$. Since $c^m=1$ on the one hand,
while $\lam$ is non-periodic on the other hand, $c\in\{-1,1\}$ as
claimed. Finally, since $-z\in\R$ whenever $z\in\R$, we conclude that
$r\le 2$ and $F(z)\subset\{-z,z\}$ for every $z\in\R$, by virtue of
analytic continuation.\qed

\bigskip
\noindent{\bf 3.7 Proposition.} The kernel of $H_\gam(z)$is
one-dimensional for every $z\in\tr\bsl\Om$.

\bigskip
\Proof: For every $z\in\tr\bsl\Om$ choose $t_n(z)$ as in Proposition
3.2, that is,
$$H_\gam(z)\eta=0~\hbox{if and only if}~k_\gam D_{G(z)}\xi
=\Gam(z)D_{G(z)^{-1}}\xi~,$$
whenever $\eta\in\ell^2(\Z)$ and
$\xi_n=(t_n(z)+t_{n}(z)^{-1})^{-1}\eta_n$.

Let $\wth_\gam(G(z))$ be defined as follows:

$$\(\wth_\gam(G(z))\eta\)_n=s_n(z)\gam\eta_{n+1}+s_n(z)^{-1}
\gam^{-1}\eta_{n-1}+\beta\(G(z)\lam^n+G(z)^{-1}\lam^{-n}\)\eta_n~,$$
where
$$s_n(z)=\(t_{n+1}(z)+t_{n+1}(z)^{-1}\)\(t_n(z)+t_n(z)^{-1}\)^{-1}~.$$
Then $\wth_\gam(G(z))$ is a bounded operator. Since $\Gam$ is analytic,
$\wth_\gam(G(z))$ is analytic in $z$. Let
$$B(z)=Q(z)\wth_\gam(G(z))Q(z)~,\quad z\in\tr\bsl\Om~.$$
Then $B(z)$ is a finite rank operator and $B$ is analytic by Lemma 3.5.
By Corollary 3.4 the cardinality of the spectrum of $B(z)$ equals the
rank of $B(z)$. The rank of $B(z)$, however, equals the dimension of the
kernel of $H_\gam(z)$. So, in order to settle the claim it has to be
shown that the spectrum of $B(z)$ contains exactly one point. Since
$Sp(B(\om(z))=Sp(B(z))$, we can define a multi-function $F$ on $\R\bsl
p(\Om)$,
$$F(z)=Sp(B(\tz))~,\quad p(\tz)=z~.$$
Then $F$ enjoys the properties stated in Lemma 3.6. It follows that
$F(z)\subset\{-z,z\}$ for every $z\in\R\bsl p(\Om)$. Since
$F(p(\tz))\subset Sp(h_\gam(G(\tz))$ for every $\tz\in\tr\bsl\Om$, and
since $-z$ cannot be an eigenvalue of $h_\gam(G(\tz))$ if $z$ is one, we
conclude that $Sp(B(\tz))$ contains exactly one element, as claimed.
(Compare this with the remark following Proposition 3.2.)\qed

As an immediate consequence of Proposition 3.7, we obtain the following
corollary, announced in the introductory remarks to this paragraph:

\bigskip
\noindent{\bf 3.8 Corollary.} If $z_1.z_2\in\tr\bsl\Om$ have the
property that $(G(z_1), \Gam(z_1))=(G(z_2), \Gam(z_2))$, then $z_1=z_2$.

\bigskip

\item{\bf 4.}{\bf The case $|\del|=1$}
\medskip
In this final paragraph we will describe how certain features from the
preceding two paragraphs carry over to the operator $h=h(1)$.

First, $h$ is a fixed point of the automorphism $\rhob$, and therefore
the operators $h$ and $k$ commute. This information can be used to show
the following:
$$\hbox{The spectrum of $h_\gam=h_\gam(1)$ is constant for
$\beta^{-1}\le|\gam|\le\beta$}~.\leqno{(4.1)}$$

\noindent Moreover, if $\xi$ is a solution of
$$h(x)\xi=\chi\xi~;~~\chi\in Sp(h)~,~~x\in\T~,$$
and $|\xi_n|$ grows moderately as $|n|\to\infty$, then
$$D_x\,k\,D_x\,\xi=c\,\xi$$
for some complex number $c$. Hence
$$D_{\x}\,k^{-1}\,D_{\x}\oxi=\ovl{D_x\,k\,D_x\,\xi}=\oc\,\oxi~,$$
or
$$D_x\,k\,D_x\,\oxi=\oc^{-1}\,\oxi~.$$
So, if $\xi$ is real, then $|c|=1$. This information, in conjunction with
(4.1), can be used to show the following:
$$Sp(D_x\,k_\gam\,D_x)=\T~\hbox{for $\be^{-1}\le|\gam|\le\beta$ and
$x\in\T$}~.\leqno(4.2)$$

In analogy to the operators $H_\gam(z)$ considered in paragraph 3, we
define for $\gam\in\r$ close to $1$ and $\theta,\nu\in\r$;
$$\(H_\gam(\theta,\nu)\eta\)_n=\sum^{\infty}_{j=-\infty}a_j\eta_{n+j}
+\tan\pi(\al n^2+2\theta n+\nu)\eta_n~,$$ where $|\eta_n|$ grows
moderately as $|n|\to\infty$, and
$$\tan g(\gam u)=\sum^{\infty}_{j=-\infty}a_ju^j~,$$
with $g$ as in (1.6). Then $H_\gam(\theta,\nu)$ determines an unbounded
operator on $\ell^2(\Z)$ with a dense domain, as long as
$$(\theta,\nu)\notin\Om_0=\left\{(\tth,\tnu)~\big|~\al n^2
+2\tth n+\tnu\in{1\over 2}\Z\right\}~.$$ As $\Gam(z)$ and $G(z)$ approach
suitable complex numbers of modulus one, respectively, $H_\gam(z)$ is
seen to approach an operator of this type, on the linear subspace of all
vectors with finitely many non-vanishing components only, say.

As in Proposition 3.2, one can now show the following: Let
$\xi\in\ell^2(\Z)$ be an eigenvector of $h_\gam(x)$ and suppose that
$$D_x\,k_\gam\,D_x\,\xi=c\xi~.$$
Let
$$\eta_n=\cos\pi(\al n^2+2\theta n+\nu)\xi_n~,$$
where $e^{2\pi\theta i}=x~,~~e^{2\pi\nu i}=c$. Then
$$H_\gam(\theta,\nu)\eta=0~,$$
in case $(\theta,\nu)\notin\Om_0$. Thus, information about the eigenvalue
problem for the operator $h_\gam(\theta)$ can be transformed into
information regarding the kernel of an operator of the form
$H_\gam(\theta,\nu)$, and to some degree, this works in the opposite
direction as well.

If $\beta$ is sufficiently large, then we may choose $\gam=1$. In this
case, $H(\theta,\nu)=H_1(\theta,\nu)$ becomes an essentially self-adjoint
operator. A complete analysis of the spectral properties of operators of
this type can be obtained in case $\al=0$ and $\theta$ is an irrational
number satisfying a stronger diophantine condition than the one that has
been assumed to be valid for $\al$ in this paper (cf. [PF], Chapter VII,
\S18).

\vfill\eject

\noindent{\bf Appendix}
\medskip
Two items related to a system of difference equations which have been
used explicitly in the text will be assembled. For more details, we refer
to [R1] and [R2].

\item{\bf A1.} Consider the system

$$\left\{\eqalign{\cos(\pi\al q+\theta)(X_{p-1,q}&+X_{p+1,q})
+\beta\cos(\pi\al p+\theta)(X_{p,q-1}
+X_{p,q+1})=zX_{pq}\cr
\sin(\pi\al q+\theta)(X_{p-1,q}&-X_{p+1,q})-\beta\sin(\pi\al p+\theta)
(X_{p,q-1}-X_{p,q+1})=0~.}\right.\leqno(A1.1)$$ As in [R1], one can
construct a recursion
$$\pmatrix{X_{p+2,p+2}\cr X_{p+2,p+1}\cr X_{p+1,p+1}}=F_p\pmatrix{X_{p+1,p+1}\cr
X_{p+1,p}\cr X_{pp}}~,
\quad p\ge 0~,$$
where $\{X_{pq}\}$ is any solution of (A1.1) and
$$F_p=E_pD_pC_p~,$$
$$\eqalign{C_p&=\left[\matrix{-\displaystyle{\beta\sin[2\pi\al(p+1)+2\theta]
\over{\sin\pi\al}}
&\displaystyle{\chi\sin[\pi\al(p+1)+\theta]\over{\sin\pi\al}}
&-\displaystyle{\sin[\pi\al(2p+1)+2\theta]\over{\sin\pi\al}}\cr
\noalign{\vskip10pt} 1&0&0\cr
\noalign{\vskip10pt}
0&1&0}\right]\cr
\noalign{\vskip15pt}
D_p&=\left[\matrix{0&\displaystyle{\chi\over{2\cos[\pi\al(p+1)+\theta]}}
&-\beta\cr
\noalign{\vskip10pt}
0&1&0\cr
\noalign{\vskip10pt}
1&0&0}\right]\cr
\noalign{\vskip15pt}
E_p&=\left[\matrix{-\displaystyle{\chi\sin[\pi\al(p+1)+\theta]
\over{\beta\sin[\pi\al(2p+3)+2\theta]}}
&-\displaystyle{\sin[2\pi\al(p+1)+2\theta]\over{\beta\sin[\pi\al(2p+3)+2\theta]}}
&\displaystyle{\sin\pi\al\over{\sin[\pi\al(2p+3)+2\theta]}}\cr
\noalign{\vskip10pt} 1&0&0\cr
\noalign{\vskip10pt}
0&1&0}\right]~.}$$ If $\{X_{pq}\}$ and $\{Y_{pq}\}$ are solutions of
(A1.1), then
$$X_{p+1,p+1}Y_{p+2,p+2}-X_{p+2,p+2}Y_{p+1,p+1}=
-\beta\displaystyle{\sin(\pi\al+\theta)\over{\sin[\pi\al(2p+3)+2\theta]}}
\cdot\Bigl(X_{00}Y_{11}-X_{11}Y_{00}\Bigr)~.$$

Now suppose $\theta\notin\pi\al\Z$. Then the following holds true: If
$(X_{00},X_{11})\not=(0,0)$, $(Y_{00},Y_{11})\not=(0,0)$ and the sequence
$$\{X_{pp}Y_{p+1,p+1}-X_{p+1,p+1}Y_{pp}\}_{p\ge 0}$$
is bounded, while $\{\sup\limits_{p\in{\bs Z}}|X_{pq}|\}_{q\ge 0}$ and
$\{\sup\limits_{p\in{\bs Z}}|Y_{pq}|\}_{q\le 0}$ are bounded as well,
then $\{X_{pq}\}$ and $\{Y_{pq}\}$ must be linearly dependent.

\bigskip

\item{\bf A2.} Now consider the system (A1.1) for $\theta=0$. Let
$$(h-z)^{-1}=\sum_{p,q\in{\bs Z}}c_{pq}(z)w_{pq}~,\quad z\in\R=\C\bsl Sp(h)~.$$
Then $c_{pq}$ is analytic and $\{c_{pq}(z)\}$ solves the system (A1.1),
except for $p=q=0$.

By [R2], Paragraph 4, there exist polynomials $d_{pq}$ which are either
zero or of degree $\bigl||p|-|q|\bigr|-1$ such that $\{d_{pq}(z)\}$
solves the system (A1.1) except for $p=q=0$ and
$$d_{pq}(z)=0~\hbox{for}~q\ge-|p|;~~~ d_{p,-p-1}(z)=(-1)^p
\beta^{-p-1}~\hbox{for}~p\ge 0~.$$
If $\{X_{pq}\}$ is a solution of (A1.1) for some $z\in\R$ which decays
exponentially uniformly in $p$ as $q\to\infty$, then $\{X_{pq}\}$ and
$\{c_{pq}(z)-d_{pq}(z)\}$ must be linearly independent.

\vfill\eject

\normalbaselines

\bigskip
\noindent{\bf References}
\medskip
\item{[AS]} L. V. Ahlfors, L. Sario, ``Riemann Surfaces'', Princeton
Mathematical Series, 1976.
\medskip
\item{[Ko]} T. Kato, ``Perturbation Theory for Linear Operators'' Second
Edition, Springer-Verlag 1976.
\medskip
\item{[Kn]} Y. Katznelson, ``An Introduction to Harmonic Analysis'',
Second Corrected Edition, Dover 1976.
\medskip
\item{[PF]} L. Pastur, A. Figotin, ``Spectra of Random and
Almost-Periodic Operators'', Springer-Verlag 1992.
\medskip
\item{[R1]} N. Riedel, ``Almost Mathieu operators and rotation
$C^*$--algebras'', Proc. London Math. Soc. 3 (56), (1988), 281--302.
\medskip
\item{[R2]} N. Riedel, ``The spectrum of a class of almost periodic operators'',
submitted since July 1992.
\medskip
\item{[R3]} N. Riedel, ``Regularity of the spectrum for the almost
Mathieu operator'' (September 1992), Proc. Amer. Math. Soc., to appear.

\bigskip

Department of Mathematics

Tulane University

e-mail:~~{\it nriedel@mailhost.tcs.tulane.edu}

\bye